\documentclass[12pt]{article}
\usepackage[utf8]{inputenc}
\usepackage[english]{babel}
\usepackage{graphicx}
\usepackage{amssymb}
\usepackage{amsthm}
\usepackage{amsmath}

\usepackage{subcaption}
\usepackage[font=footnotesize]{caption}

\usepackage{mathtools}

\usepackage{tikz}
\usetikzlibrary{matrix}
\usepackage[left=1.5in,right=1in]{geometry}

\usepackage[T1]{fontenc}

\usepackage[shadow, colorinlistoftodos,textsize=tiny]{todonotes}
\setuptodonotes{fancyline, backgroundcolor=gray!10,bordercolor=gray}

\usepackage{authblk}

\usepackage[capitalise,noabbrev, nameinlink, poorman]{cleveref}

\newtheorem{teo}{Theorem}[section]

\newtheorem{coro}[teo]{Corollary}
\newtheorem{prop}[teo]{Proposition}

\newtheorem{lema}[teo]{Lemma}
\theoremstyle{definition}

\newtheorem{ej}[teo]{Example}

\newtheorem{preg}{Question}

\newenvironment{dem}[1][Proof]{\begin{trivlist}
  \item[\hskip \labelsep {\bfseries #1}]}{\end{trivlist}}

\theoremstyle{remark}

\newcommand{\N}{\mathbb{N}}

\newcommand{\floor}[1]{\left \lfloor #1 \right \rfloor}

\newcommand{\cP}{\mathcal{P}}

\newcommand{\cO}{\mathcal{O}}

\newcommand{\cM}{\mathcal{M}}
\newcommand{\cT}{\mathcal{T}}

\newcommand{\cK}{\mathcal{K}}
\newcommand{\cX}{\mathcal{X}}
\newcommand{\cY}{\mathcal{Y}}
\newcommand{\cZ}{\mathcal{Z}}

\newcommand{\ol}[1]{\overline{#1}}
\newcommand{\os}[1]{\widetilde{#1}}

\newcommand{\gen}[1]{\langle #1 \rangle} 

\newcommand{\Mon}{{\rm Mon}}
\newcommand{\1}{{\mathbf{1}}}
\newcommand{\2}{{\mathbf{2}}}

\newcommand{\ertimes}{\rtimes_\eta}
\newcommand{\XY}{\cX\ertimes\cY}
\newcommand{\MY}{\cM\ertimes\cY}
\newcommand{\opp}[1]{#1^{\rm opp}}

\title{Semiregular abstract polyhedra with trivial facet stabilizer
}
\author{Elías Mochán}
\affil{Department of Mathematics, Northeastern University, 02115
  Boston, USA\thanks{\tt email:jaime.mochan@im.unam.mx}}
\begin{document}
\maketitle


\abstract{\emph{Abstract polytopes} generalize the face lattice of
  convex polytopes. A polytope is \emph{semiregular} if
  its facets are regular and its automorphism group acts transitively
  on its vertices. In this paper we construct semiregular,
  facet-transitive polyhedra with trivial facet stabilizer, showing
  that semiregular abstract polyhedra can have an unbounded number of
  flag orbits, while having as little as one facet orbit. We interpret
  this construction in terms of operations applied to high rank
  regular and chiral polytopes, and we see how these same operations
  help us construct alternating semiregular polyhedra (that is, with
  two facet orbits and adjacent facets in different orbits). Finally,
  we give an idea to generalize this construction giving examples in
  higher ranks.}

\section{Introduction}\label{sec:Intro}

Abstract polytopes are combinatorial objects that generalize the face
lattice of convex polytopes and other geometric objects. Just as with
geometric polytopes, one of the main approaches to study abstract
polytopes has been through their symmetry. A \emph{flag} is typically
defined as a set having one face of each rank, and satisfying that all
those faces are incident to each other. One can prove that the
automorphism group of an abstract polytope acts semiregularly on its
flags, and therefore, one can measure how symmetric a polytope is by
counting its flag orbits: the fewer flag orbits the more
symmetrical. Regular polytopes are those with exactly one flag orbit,
and are by far the most studied class of abstract polytope, having a
whole book~\cite{ARP} dedicated to them.

When studying other classes of polytopes, some have studied polytopes
with a small number of flag orbits (see for example~\cite{3Orb},
\cite{2OrbYo}, \cite{2OrbMani}, \cite{QuiralesEgonAsia}), but some
classes where the number of flag orbits is not restricted have been
studied too. Such is the case of semiregular polytopes: those with
regular facets and where the automorphism group acts transitively on
the vertices. Despite the number of flag orbits being unbounded
(see~\cite{OrbitasEnSemiregs}), most of the research on semiregular
polytopes has been focused on constructions where the result has only
one or two flag orbits (see~\cite{SemiReg}, for example). This is due
to the fact that the examples constructed are hereditary: every
symmetry of a facet induces a symmetry of the whole polytope.

In this paper we construct semiregular polyhedra that are
face-transitive and that are the opposite of hereditary: no
non-trivial symmetry of a facet induces a symmetry of the whole
polyhedron. In other words, our examples have trivial facet
stabilizer. This implies that despite having one orbit on vertices and
one orbit on faces, the number of orbits on flags in unbounded in our
construction.

For our construction it is convenient to think of polytopes as a
special class of maniplexes: graphs with a proper coloring of their
edges satisfying certain conditions. Thinking of polytopes in terms of
graphs allows us to use voltage graphs to construct them, and thanks
to the results of~\cite{CayExt} and~\cite{IntPropYo}, we know
conditions on the voltage assignments that tell us when the graphs we
obtain represent a polytope.

We will also see that our examples can be interpreted as the result of
applying some particular operations to regular and chiral polytopes
(those with the highest possible rotational symmetry but no
reflections) of high rank. These same operations can be used to obtain
alternating semiregular polyhedra with trivial facet stabilizer (that
is, semiregular polyhedra with two face-orbits and where each edge
lies between faces in different orbits).

In \cref{sec:basic} we introduce the basic notions necessary for this
paper. \cref{sec:SymType} presents the concept of \emph{symmetry type
  graph}, which is used later in the paper. \cref{sec:volts}
presents the non-standard definition of \emph{graph} used in this
paper, together with the concept of voltage graphs and how to
construct the derived graph of a voltage graph. In \cref{sec:SymTypes}
we characterize the regular symmetry type graphs that semiregular
facet-transitive polyhedra with trivial facet stabilizer may
have. \cref{sec:fam1const} and \cref{sec:fam2const} deal with how to construct
polyhedra with the symmetry types determined in
\cref{sec:SymTypes}. In \cref{sec:VoltOps} we interpret our
construction in terms of operations. In \cref{sec:AltSemiReg} we see
how to use these operations to get examples of alternating semiregular
polyhedra. Finally, in \cref{sec:highrank} we see how we could
generalize some of our constructions to get semiregular,
facet-transitive maniplexes with trivial facet stabilizer in higher
rank, although we are not sure yet whether they are
polytopes or not.

\section{Basic concepts}\label{sec:basic}

Abstract polytopes, introduced by Schulte and Danzer
in~\cite{EgonPhD}, are typically defined as posets with certain
properties that make them similar to the face lattice of a convex
polytope. However, thanks to~\cite[Theorem 5.3]{PolyMani}, we know
that they can be equivalently defined as a particular kind of
maniplex. For simplicity, we will use this other definition.

An \emph{$n$-maniplex} (or \emph{maniplex of rank $n$}) is a connected
simple graph $\cM$ with a coloring of its edges with colors in
$[0,n-1]:=\{0,1,\ldots,n-1\}$ such that

\begin{itemize}
\item every vertex is incident to exactly one edge of each color, and
\item given two colors $i$ and $j$ such that $|i-j|>1$, the
  alternating paths of length 4 using edges of these two colors are
  closed.
\end{itemize}
The vertices of a maniplex are referred to as \emph{flags}. For the
purposes of this paper all maniplexes are finite, that is, they
have a finite number of flags. Maniplexes are a natural generalization
of maps to higher dimensions. To see other definitions of maniplex
and also learn the topological interpretation of this definition
see~\cite{Maniplexes}.

Given an $n$-maniplex $\cM$ and a color $i\in [0,n-1]$, we denote by
$\cM_{\ol{i}}$ the subgraph of $\cM$ obtained by removing the edges of
color $i$. The connected components of $\cM_{\ol{i}}$ are the
\emph{$i$-faces of $\cM$}. The $i$-face containing the flag $\Phi$ is
denoted by $(\Phi)_i$, and the set of all $i$-faces of $\cM$ is
denoted by $(\cM)_i$. The $(n-1)$-faces of an $n$-maniplex are its
\emph{facets}, while the 0-faces and 1-faces are its \emph{vertices}
and \emph{edges} respectively.


We say that an $i$-face $F$ and a $j$-face $G$ of a maniplex are
\emph{incident} if they have non-empty intersection. If in addition
$i\leq j$, we write $F\leq G$. In this way the faces of a
maniplex form a poset (see~\cite[Proposition~3.1]{PolyMani}).

\begin{figure}[h!]
  \centering
  \includegraphics[width=0.5\linewidth]{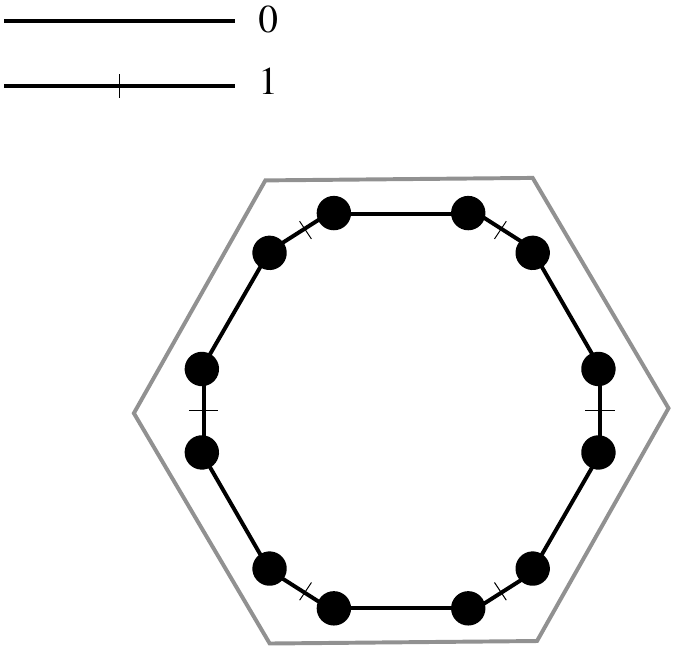}
  \caption{A hexagon viewed as a 2-maniplex.}
  \label{fig:Hexagon}
\end{figure}

\begin{figure}[h!]
  \centering
  \includegraphics[width=0.6\linewidth]{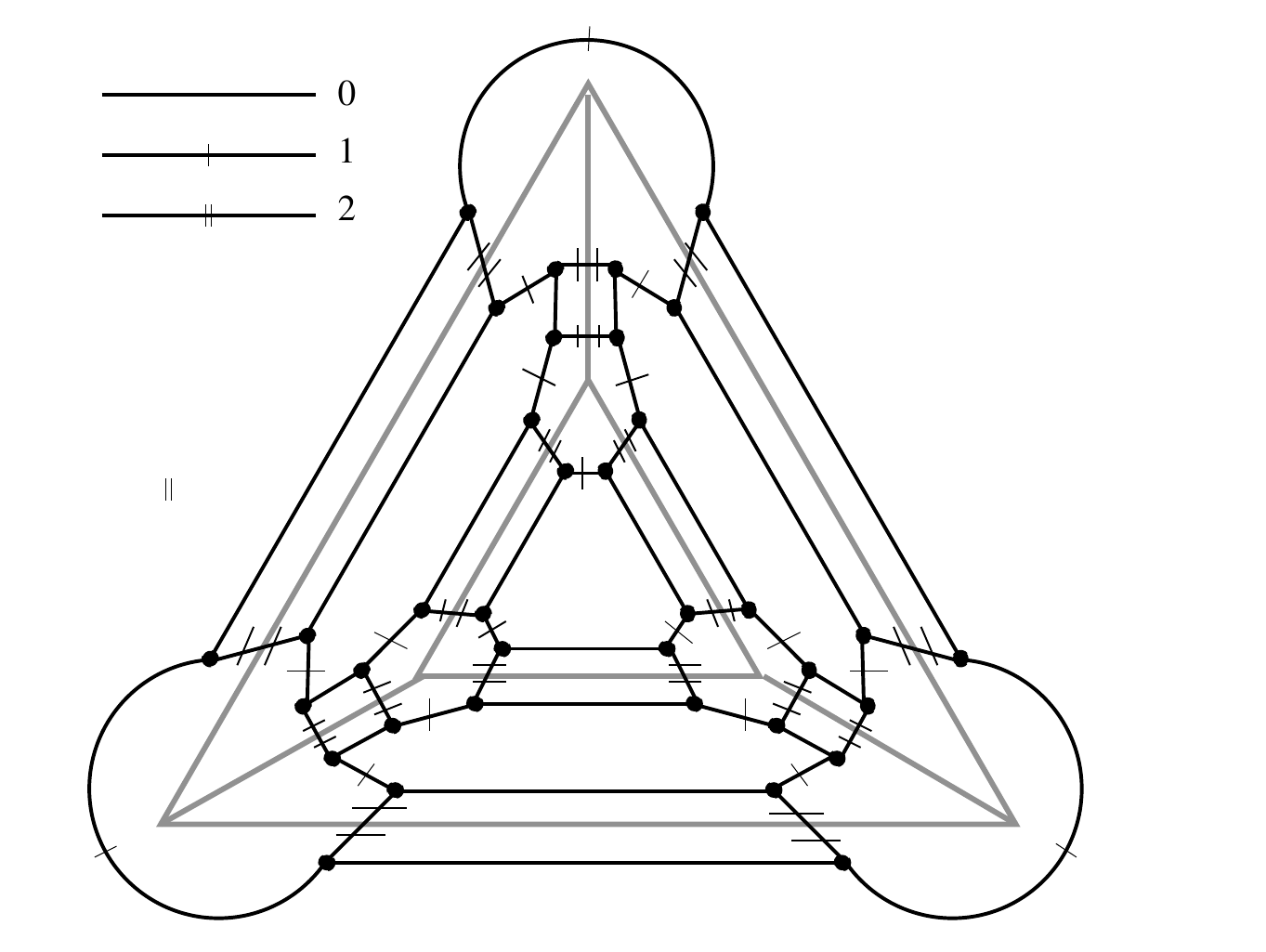}
  \caption{A prism viewed as a 3-maniplex.}
  \label{fig:Prism}
\end{figure}

In \cref{fig:Hexagon} we see how a hexagon can be interpreted as a
2-maniplex. Note that the vertices of the hexagon correspond to the
edges of color 1 of the graph, and the edges of the hexagon correspond
to the edges of color 0. Analogously, in \cref{fig:Prism} we see a
triangular prism interpreted as a 3-maniplex. Here the vertices of the
prism correspond to cycles with edges of colors 1 and 2 in the graph,
the edges of the prism correspond to cycles with edges of colors 0 and
2 in the graph, and the 2-faces of the prism correspond to cycles
with edges of colors 0 and 1 in the graph.

Given a maniplex $\cM$ we define its \emph{dual} $\cM ^*$ as the
maniplex we get by replacing the edges of color $i$ with edges of
color $n-1-i$ for each $i\in[0,n-1]$. Note that the vertices of
$\cM^*$ correspond to the facets of $\cM$ and vice-versa. If $\cM^*$
is isomorphic to $\cM$ we say that $\cM$ is \emph{self-dual}.

An \emph{abstract $n$-polytope}, or \emph{abstract polytope of rank
  $n$}, is an $n$-maniplex $\cP$ satisfying the following condition:

{\bf Path intersection property:} Let $\Phi$ and $\Psi$ be two flags of
$\cP$. If there is a path from $\Phi$ to $\Psi$ that uses only colors
in $[0,m]$ and another path from $\Phi$ to $\Psi$ that uses only
colors in $[k,n-1]$, then there is a path from $\Phi$ to $\Psi$ that
uses only colors in $[k,m]$.

In~\cite[Theorem 5.3]{PolyMani} it is proved that this definition is
equivalent to the poset definition given by Shculte and Danzer in
~\cite{EgonPhD}. From here on, whenever we talk about polytopes in this
paper we will be thinking about abstract polytopes.

An important property of polytopes when thought of as maniplex is that
the flag $\Phi^i$, $i$-adjacent to the flag $\Phi$ has the same
$j$-face as $\Phi$ for every $j\neq i$, but it has a different
$i$-face. In fact, in the study of abstract polytopes this is usually
the definition of $i$-adjacent flags.


An \emph{isomorphism} of $n$-maniplexes is just a graph isomorphism
that preserves the colors of the edges. Naturally, an automorphism of
a maniplex is just an isomorphism onto itself. The automorphism group
of a maniplex $\cM$ is denoted by $\Gamma(\cM)$. It is known that the
automorphism group of a maniplex acts freely on its flags (the proof
is identical to that of~\cite[Proposition~2A4]{ARP}, which is the
particular case of abstract polytopes). We say that a maniplex $\cM$
is \emph{regular} (some authors prefer \emph{reflexible}) if
$\Gamma(\cM)$ acts transitively (and therefore regularly) on the flags
of $\cM$.

A maniplex is called \emph{semiregular} if it is vertex-transitive
(that is, transitive on its 0-faces) and each of its facets is a
regular $(n-1)$-maniplex. Note that, since all polygons (2-polytopes)
are regular, for polyhedra (3-polytopes) semiregularity is the same as
vertex-transitivity.




Given a maniplex $\cM$, we define $r_i$ as the flag permutation that
maps each flag to its $i$-adjacent flag. It is easy to see that $r_i$
is an involution and that, if $|i-j|>1$, then $r_ir_j$ is also an
involution. In particular $r_i$ and $r_j$ commute whenever $i$ and $j$
are not consecutive. If there is an $i\in [0,n-2]$ such that
$(r_ir_{i+1})^2$ fixes some flag, we say that $\cM$ is \emph{degenerate},
otherwise it is \emph{non-degenerate}.

The group $\Mon(\cM)=\gen{r_i:0\leq i \leq n-1}$ is called the
\emph{monodromy group of $\cM$}, and we will call each of its elements
a \emph{monodromy}. It is worth mentioning that some authors prefer
the term \emph{connection group} and the notation
$\mathrm{Con}(\cM)$. An alternative way of defining the automorphism
group of a maniplex $\cM$ is as the set of flag permutations that
commute with the elements of the monodromy group.

For this paper, the automorphism group will act always on the
right-hand side and the monodromy group on the left-hand side.

If $\cM$ is a regular maniplex,
for a fixed flag $\Phi$ we can define $\rho_i$ as the (unique)
automorphism that maps $\Phi$ to $\Phi^i$. It is known that
$\{\rho_i:i\in\{0,1,\ldots,n-1\}\}$ is a generating set for the
automorphism group of $\cM$ (see~\cite[Theorem 2B8]{ARP} for the
polytopal case). Moreover, the function mapping $\rho_i$ to $r_i$ can
be extended to a group isomorphism between $\Gamma(\cM)$ and
$\Mon(\cM)$ (see, for example,~\cite[Theorem~3.9]{MixAndMon}
or~\cite[Section~7]{Maniplexes}). In particular, the monodromy group
of a regular maniplex acts regularly on its flags. It is also worth
noticing that a regular maniplex is degenerate if and only if
$\rho_i\rho_{i+1}$ is an involution for some $i\in[0,n-2]$.

\section{Symmetry type graphs}\label{sec:SymType}

For this paper we need to think of graphs in a way that is not
completely standard. A \emph{graph} $\cX$ consists of a set of
\emph{vertices} $V(\cX)$, a set of \emph{darts} $D(\cX)$, a function
$I:D(\cX)\to V(\cX)$ and an involution
$(\cdot)^{-1}:D(\cX)\to D(\cX)$. We call $I(d)$ the \emph{initial
  vertex of $d$} and $d^{-1}$ the \emph{inverse} of $d$. An
\emph{edge} of $\cX$ is just a set of the form $\{d,d^{-1}\}$.  A
\emph{semi-edge} is a dart $d$ satisfying $d=d^{-1}$. If $d$ is a
semi-edge, we abuse language and call the edge $\{d\}$ a semi-edge as
well. We think of the edge $e=\{d,d^{-1}\}$ as a line joining the
vertices $I(d)$ and $T(d):=I(d^{-1})$ (called the \emph{final vertex
  of $d$}), and we think of its darts as the two (or one) possible
orientations. A \emph{path} $W$ on $\cX$ is a finite sequence
$d_1d_2\cdots d_\ell$ of darts in $\cX$ such that $T(d_i)=I(d_{i+1})$
for all $i\in \{0,1,\ldots,\ell-1\}$. The vertex $I(d_1)$ is the
\emph{initial vertex of $W$} while $T(d_\ell)$ is its \emph{final
  vertex}. The number $\ell$ is called the \emph{length of $W$}. We
refer the reader to~\cite{Voltajes2} for details.

The \emph{symmetry type graph} of a polytope (or maniplex) $\cP$,
denoted by $\cT(\cP)$, is a graph whose vertices are the flag-orbits
of $\cP$ and satisfying that it has an edge of color $i$ between the
orbit of $\Phi$ and the orbit of $\Psi$ if and only if $\Psi$ is in
the same orbit as $\Phi^i$. If $\Phi$ and $\Phi^i$ are in the same
orbit, we draw a semi-edge of color $i$ at the vertex representing
that orbit in the symmetry type graph.

The symmetry type graph of an $n$-polytope or $n$-maniplex is an
example of an \emph{$n$-premaniplex}, defined as a graph with a
coloring of its edges with colors in $\{0,1,\ldots,n-1\}$ such that
each vertex is incident to exactly one edge of each color, and
satisfying that whenever $|i-j|>1$, the alternating paths of length 4
with colors $i$ and $j$ are closed. For the purposes of this paper all
premaniplexes are considered to be connected. In contrast to a
maniplex, a premaniplex need not be simple. The notions of
isomorphism, automorphism, regularity, monodromy, and so on, extend
naturally from maniplexes to premaniplexes.

Given two $n$-premaniplexes $\cM$ and $\cX$, a
\emph{covering projection} (\emph{covering} for short) from $\cM$ to
$\cX$, is a function $p$ from the vertices of $\cM$ to the vertices of
$\cX$ that preserves $i$-adjacencies. It can be proved that all
coverings are surjective. If $e$ is an edge of $\cM$ connecting vertices
$x$ and $y$, we define $p(e)$ as the edge connecting $p(x)$ and $p(y)$
that has the same color as $e$. If there is a covering from $\cM$ to
$\cX$, we say that \emph{$\cM$ covers $\cX$}. Note that if $\cM$ is a
maniplex and $\cX$ is its symmetry type graph, the natural projection
from $\cM$ to $\cX$ is always a covering.

Some notable symmetry type graphs are the following:
\begin{itemize}
\item The symmetry type graph of regular $n$-polytopes is denoted by
  $\1^n$. It has exactly one vertex and semi-edges of each color from
  0 to $n-1$ at that vertex.
\item The symmetry type graph of chiral $n$-polytopes is denoted by
  $\2^n_\emptyset$: It has exactly two vertices and they are joined by
  one edge of each color from 0 to $n-1$.
\item More generally, if $\cP$ is a polytope with two flag orbits and
  $I\subset \{0,1,\ldots,n-1\}$ is the set of colors such that $\Phi$
  and $\Phi^i$ are in the same orbit, then the symmetry type graph of
  $\cP$ is $\2^n_I$, which has exactly two vertices joined by one edge
  of each color in $\ol{I}$ and with semi-edges of each color in $I$
  at each vertex (see~\cite[Section 4]{SymType}).
\end{itemize}

The symmetry type graph of a polytope can be used to know how many
orbits the polytope has on $i$-faces. In fact, if $\cX=\cT(\cP)$ then
each connected component of $\cX_{\ol{i}}$ (the graph obtained from
$\cX$ by removing the edges of color $i$), corresponds to an orbit of
$i$-faces (see~\cite[Proposition 3.2]{SymType}). In particular, $\cP$ is
$i$-face-transitive if and only if $\cX_{\ol{i}}$ is connected.

A polytope is \emph{hereditary} if every symmetry of each facet
induces a symmetry of the whole polytope, or in other words, if the
stabilizer of each facet coincides with its automorphism group.  In
terms of the symmetry type graph, a polytope $\cP$ with symmetry type
graph $\cX$ is hereditary if each connected component of
$\cX_{\ol{n-1}}$ is isomorphic to the symmetry type graph of a facet
in the corresponding orbit. In particular, if $\cP$ has regular facets
and is hereditary, then each connected component of $\cX_{\ol{n-1}}$
must have only one vertex, which in turn means that $\cX$ is either
$\1^n$ (implying that $\cP$ is regular), or $\2^n_{\ol{n-1}}$: the
graph with two vertices joined by an edge of color $n-1$ and
semi-edges of all other colors at each vertex (for more on this,
see~\cite{Hereditary}). Note that if we remove edges of any color
other than $n-1$ from these two graphs we do not disconnect them,
meaning that hereditary polytopes with regular facets are $i$-face
transitive for all $i<n-1$.  Polytopes with symmetry type graph
$\2^n_{\ol{n-1}}$ are a particular case of \emph{alternating
  semiregular polytopes}: semiregular polytopes with two facet orbits
where $(n-1)$-adjacent flags are in different orbits. This kind of
polytope has been studied in~\cite{SemiReg}, for example.

In this paper we want to focus on semiregular polytopes that are the
opposite of hereditary: Instead of wanting every symmetry of a facet
to induce a symmetry of the whole polytope, we want that none of them
do; that is, we want the stabilizer of every facet to be trivial. This
means that any two flags in the same facet must be in different
orbits. In terms of the symmetry type graph, this means that if $\cP$
is a polytope with the desired property and $\cX$ is its symmetry type
graph, each connected component of $\cX_{\ol{n-1}}$ must be isomorphic
to a facet in the corresponding orbit.

Note that if a hereditary polytope with regular facets is not regular,
it must have exactly two facet orbits. This is not the case when we
remove heredity: semiregular polytopes can have any number of facet
orbits, even if all facets are isomorphic to each other
(see~\cite[Theorem 5.3]{OrbitasEnSemiregs}).

In~\cite{OrbitasEnSemiregs}, Pisanski, Schulte and Ivić Weiss
construct examples of semi-regular $n$-polytopes with trivial vertex
stabilizer and show that they can have arbitrarily many flag
orbits. However, their examples have (usually) many facet-orbits (even
though they are isomorphic to each other). We will show in this
article that it is also possible to construct facet-transitive
semiregular polytopes that are not only not regular, but that have as
many flag-orbits as the number of flags in a facet.


\section{Voltage graphs}\label{sec:volts}

Voltage graphs are a tool used in graph theory to construct graphs
with a given quotient (see~\cite{VoltsLibro}). They have recently been
used to construct maniplexes and polytopes with some prescribed
symmetry types (see~\cite{3Orb} and~\cite{2OrbYo}), as we will do in
this paper.

A \emph{voltage graph} is a triplet $(\cX,\Gamma,\xi)$ where $\cX$ is
a graph, $\Gamma$ is a group (called the \emph{voltage group}), and
$\xi:D(\cX)\to\Gamma$ is a function (called \emph{voltage
  assignment}), satisfying that $\xi(d^{-1})=\xi(d)^{-1}$ for every
dart $d\in D(\cX)$. The element $\xi(d)$ is called the \emph{voltage
  of $d$}. Moreover, if $W=d_1 d_2 d_3\cdots d_k$ is a path in $\cX$,
then we define $\xi(W)$ as $\xi(d_k)\cdots\xi(d_2)\xi(d_1)$, and call
this the \emph{voltage of $W$}. Note that if $W$ ends at the initial
vertex of $W'$, then $\xi(WW')=\xi(W')\xi(W)$.

Given a voltage graph $(\cX,\Gamma,\xi)$ we can construct a graph
$\cX^\xi$ where $\Gamma$ acts by automorphisms. This action is
semiregular on both the vertices and darts of $\cX^\xi$ and the
quotient $\cX^\xi/\Gamma$ is isomorphic to $\cX$. To define $\cX^\xi$
we make $V(\cX^\xi)=V(\cX)\times \Gamma$ and
$D(\cX^\xi)=D(\cX)\times \Gamma$. Then we define the initial vertex of
the dart $(d,\gamma)$ to be $I(d,\gamma)=(I(d),\gamma)$; and the
inverse of $(d,\gamma)$ to be
$(d,\gamma)^{-1} = (d^{-1},\xi(d)\gamma)$. In other words, if there is
a dart $d$ from $x$ to $y$ in $\cX$, then for each $\gamma\in\Gamma$
there is a dart $(d,\gamma)$ from $(x,\gamma)$ to $(y,\xi(d)\gamma)$
in $\cX^\xi$.  If $\cX$ has a coloring (of its vertices, darts
and\slash or edges) we color $(x,\gamma)$ with the same color as $x$.

One can verify that the $\Gamma$-action
$(x,\gamma)\tau = (x,\gamma\tau)$ (where $x$ might be either a vertex
or a dart) on $\cX^\xi$ is by (color-preserving) automorphisms,
semiregular and that two elements (darts or vertices) are in the same
$\Gamma$-orbit if and only if they have the same first
coordinate. Therefore, $\cX^\xi/\Gamma$ is isomorphic to $\cX$. This
gives a natural embedding of $\Gamma$ into the automorphism group of
$\cX^\xi$.

It is convenient to work with voltage graphs where $\cX$ has a
spanning tree $T$ satisfying that all its darts have trivial
voltage. In this case, each connected component of $\cX^\xi$ has
$k|V(\cX)|$ vertices, where
$k=|\gen{\xi(d):d\in D(\cX)\setminus D(T)}|$. To be more precise, the
connected component of the vertex $(x,\gamma)$ has vertex set
$V(\cX)\times \gamma\gen{\xi(d):d\in D(\cX)\setminus D(T)}$.

A \emph{voltage premaniplex} is a voltage graph $(\cX,\Gamma,\xi)$
where $\cX$ is a premaniplex and $\xi$ satisfies the following
properties:

\begin{itemize}
\item $\Gamma = \gen{\xi(d):d\in D(\cX)\setminus D(T)}$ for a spanning
  tree $T$ of $\cX$ with trivial voltage on all its darts.
\item The voltage of every semiedge has order 2.
\item If $d$ and $d'$ are \emph{parallel darts} (that is, they share
  their initial and final vertices) then $\xi(d)\neq \xi(d')$.
\item Whenever $|i-j|>1$, the voltage of any alternating path of
  length 4 with colors $i$ and $j$ is trivial.
\end{itemize}

In~\cite[Lemma 3.1]{IntPropYo} it is proven that if $(\cX,\Gamma,\xi)$
is a voltage premaniplex, then $\cX^\xi$ is a maniplex. As mentioned
before, $\cX^\xi/\Gamma$ is isomorphic to $\cX$, so if one can prove
that every automorphism of $\cX^\xi$ can be interpreted as an element
of $\Gamma$, then $\cX$ is the symmetry type graph of $\cX^\xi$.

In general it can be hard to prove that $\Gamma$ is the whole
automorphism group of $\cX^\xi$. However, some automorphisms of
$\cX^\xi$ not in $\Gamma$ are easier to understand than others. An
automorphism $\os{\tau}$ of $\cX^\xi$ is a \emph{lift} of an
automorphism $\tau$ of $\cX$ if it acts as $\tau$ on the first
coordinate. In this case we also say that $\os{\tau}$ \emph{projects
  (to $\tau$)} and that $\tau$ \emph{lifts (to $\os{\tau}$)}. Note that an automorphism of $\cX$ may have many lifts
(in fact, it has either 0 or $|\Gamma|$ lifts).

The following theorem describes when an automorphism of a voltage
graph lifts.

\begin{teo}\label{teo:lifts}\cite[Theorem 4.1]{Voltajes2}
  Let $(X,\Gamma,\xi)$ be a voltage graph such that $\cX^\xi$ is
  connected. Let $x$ be a vertex of $\cX$ and $\tau$ an automorphism
  of $\cX$. Then $\tau$ 
  lifts if and only if the function that maps $\xi(W)$ to $\xi(W\tau)$
  for every closed path $W$ based at $x$ is well defined. In this
  case, this function is an automorphism of $\Gamma$.
\end{teo}

A way of interpreting the previous theorem is by saying that an
automorphism of $\cX$ lifts if and only if it \emph{induces} an
automorphism of $\Gamma$. This is very useful, as we can understand
automorphisms that lift in terms of automorphisms of the voltage
group. In particular, if we can prove that no non-trivial automorphism
of $\cX$ induces an automorphism of $\Gamma$, we will know that
the automorphisms of $\cX^\xi$ that project are precisely the elements
of $\Gamma$.

The automorphisms that do not project are, in general, much harder to
identify. Fortunately, we can avoid them completely in the case when
$\cX$ is a regular premaniplex, by using the following theorem:

\begin{teo}\label{teo:symtypereg}\cite[Corollary 5.2]{CayExt}
  Let $(\cX,\Gamma,\xi)$ be a voltage premaniplex. If $\cX$ is regular,
  then every automorphism of $\cX^\xi$ projects.
\end{teo}

\section{Regular symmetry type graphs of semiregular facet-transitive
  polyhedra with trivial facet stabilizer}\label{sec:SymTypes}

In this section we will characterize those regular premaniplexes that
are symmetry type graphs of semiregular, facet-transitive polyhedra
with trivial facet stabilizer. Recall that, since all abstract
polygons are regular, an abstract polyhedron is semiregular if and
only if it is vertex-transitive, so we are looking for polyhedra that
are both vertex- and facet-transitive.

Let us ask: What could be the symmetry type graph $\cX$ of a
facet-transitive semiregular ($n+1)$-polytope $\cP$ with trivial facet
stabilizer? As established in \cref{sec:SymType}, $\cX_{\ol{n}}$ must
be connected for $\cP$ to be facet-transitive, and it must be
isomorphic to a facet of $\cP$ to ensure that facet-stabilizers are
trivial in $\cP$. Therefore, to construct $\cX$ we may start with a
regular $n$-polytope $\cK$ (which will be isomorphic to a facet of
$\cP$), and then add edges of color $n$ in such a way that:
\begin{enumerate}
\item The resulting graph $\cX$ is an $(n+1)$-premaniplex.
\item $\cX_{\ol{0}}$ is connected (to ensure that $\cP$ is
  vertex-transitive).
\end{enumerate}

So, for $n=3$ we are looking for a 3-premaniplex $\cX$ satisfying that
$\cX_{\ol{2}}$ is a polygon, and that $\cX_{\ol{0}}$ is connected. In
order to be able to use \cref{teo:symtypereg}, we will look for
examples where $\cX$ is also regular.

\begin{figure}[h!]
  \centering
  \includegraphics[width=0.5\linewidth]{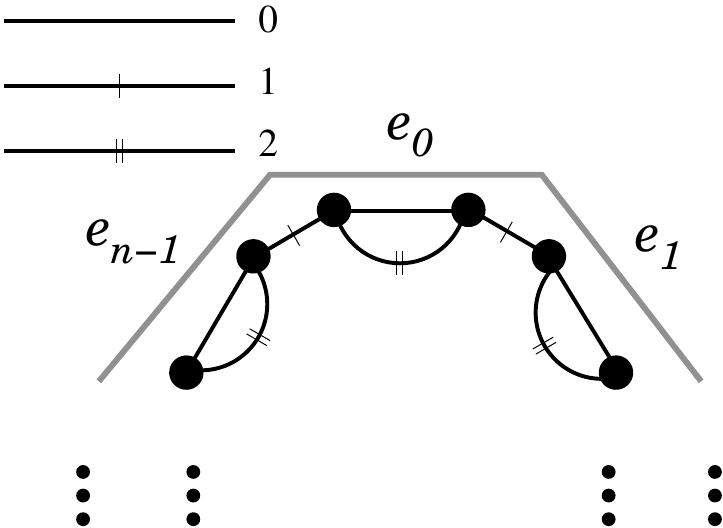}
  \caption{Symmetry type graphs of family 1.}
  \label{fig:fam1}
\end{figure}

Let us find all 3-premaniplexes with those properties:

Start with an $n$-gon, that is, a 2-maniplex with $2n$ flags forming
an alternating cycle with colors 0 and 1. We want to add edges of
color 2 that ``commute'' with the edges of color 0 in such a way that
our conditions are satisfied. First notice that if there is a
semi-edge of color 2, because of regularity, all edges of color 2
would be semi-edges. In this case $\cX_{\ol{0}}$ would have $n$
connected components (each consisting of an edge of color 1 and two
semi-edges of color 2), and this contradicts our hypothesis that
$\cX_{\ol{0}}$ is connected. Therefore, there cannot be any semi-edge
in $\cX$. Since the alternating paths of length 4 with colors 0 and 2
are closed, using regularity we get that either they are all 4-cycles
or they are all 2-cycles that are traveled twice. If these paths are
2-cycles, it means that every edge of color 2 is parallel to an edge
of color 0. If $\cX$ is an $n$-gon with additional edges of color 2
parallel to its edges of color 0, then $\cX_{\ol{0}}$ is an
alternating cycle of length $2n$ with colors 1 and 2, which is
connected. Therefore, we have found our first family of premaniplexes
with the desired properties: $n$-gons with edges of color 2 connecting
pairs of 0-adjacent flags (see \cref{fig:fam1}). We will call this
\emph{family 1}.

\begin{figure}[h!]
  \centering
  \includegraphics[width=0.5\linewidth]{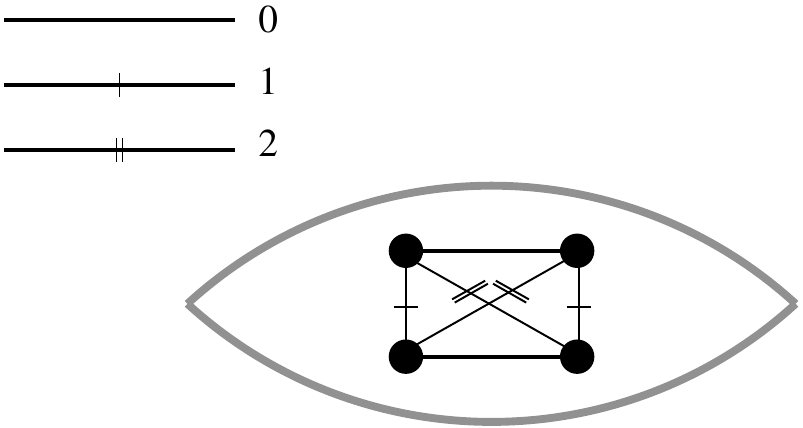}
  \caption{The sporadic premaniplex $\cX_s$.}
  \label{fig:sporadic}
\end{figure}

Now suppose that $\cX$ satisfies our conditions and it has a 4-cycle
with flags $\Phi$, $\Phi^0$, $\Phi^{02}$ and $\Phi^2$. There is an
automorphism $\tau$ of $\cX$ that maps $\Phi$ to $\Phi^0$ and $\Phi^2$
to $\Phi^{02}$. This means that $\tau$ fixes both the 1-face $e$ that
consists of $\Phi$ and $\Phi^0$ and $e'$ that consists of $\Phi^2$ and
$\Phi^{02}$. The only way two different edges of a polygon are fixed
by the same non-trivial automorphism is if they are opposites in a
polygon with an even number of sides. So we know that $n$ must be even
and that 2-adjacent flags lie on opposite edges. The edge opposite to
$e$ has the flags $(r_1r_0)^{n/2}\Phi$ and $r_0(r_1r_0)^{n/2}\Phi$,
and this edge must be $e'$. So, $\Phi^2$ must be one of these two
flags. If $\Phi^2=(r_1r_0)^{n/2}\Phi$ we get an alternating 4-cycle of
colors 1 and 2 with flags $\Phi$, $\Phi^2$, $\Phi^{21}$ and
$\Phi^1$. This means that for $n\geq 3$, $\cX_{\ol{0}}$ cannot be
connected. For $n=2$, however, we have found one sporadic example
$\cX_s$: a digon with edges of color 2 between opposite flags (see
\cref{fig:sporadic}).

Now let us investigate the case when
$\Phi^2=r_0(r_1r_0)^{n/2}\Phi$. Since we are assuming that $\cX$ is
regular we can use the fact that the monodromy group acts regularly on
the flags to get rid of $\Phi$ in our notation and just write
$r_2=r_0(r_1r_0)^{n/2}$. Then, we get that
\[
  (r_1 r_2)^2 = (r_1 r_0)^{2(n/2+1)} = (r_1 r_0)^{n+2} = (r_1 r_0)^2.
\]

Let the edges (1-faces) of $\cX$ be called $e_0, e_1, \ldots, e_{n-1}$
where $e_k$ is the edge that contains the flag
$(r_1r_0)^k\Phi$. 
We will call the flags of the
form $(r_1 r_0)^k\Phi$ \emph{left flags}, and the flags of the form
$r_0(r_1 r_0)^k\Phi$ \emph{right flags}. Notice that the alternating
path of length 4 with flags $\Phi$, $\Phi^2$, $\Phi^{21}$ and
$\Phi^{212}$ has: \
\begin{enumerate}
\item the left flag of $e_0$,
\item the right flag of $e'=e_{n/2}$,
\item the left flag of $e_{n/2+1}$, and
\item the right flag of $e_1$.
\end{enumerate}

If we continue the alternating path with colors 1 and 2 we will get

\begin{enumerate}
\item the left flags of the edges $e_{2k}$,
\item the right flags of the edges $e_{n/2+2k}$,
\item the left flags of the edges $e_{n/2+2k+1}$, and
\item the right flags of the edges $e_{2k+1}$,
\end{enumerate}
for all $k$.

\begin{figure}[h!]
  \centering
  \includegraphics[width=0.5\linewidth]{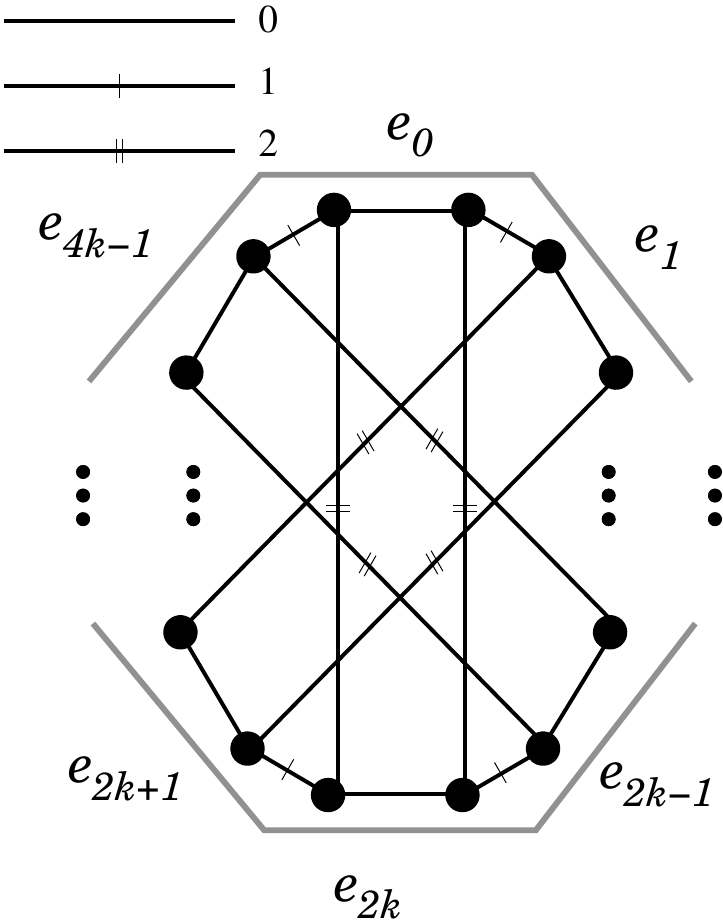}
  \caption{Symmetry type graphs of family 2.}
  \label{fig:fam2}
\end{figure}

We want this path to go through every flag before going back to
$\Phi$. The only way to make sure that all these flags are different
is if $n/2$ is even, so that every 4 steps have a left flag and a right
flag of an even edge and of an odd edge. Therefore, our second family
of premaniplexes with the desired conditions are $n$-gons where $n$ is
divisible by 4, and the edges of color 2 connect each flag to the
0-adjacent flag of its opposite flag (see \cref{fig:fam2}). We will
call this \emph{family 2}.

We have proven the following:

\begin{prop}\label{prop:regpmpx}
  Let $\cX$ be a regular 3-premaniplex such that $\cX_{\ol{2}}$ is a
  polygon and $\cX_{\ol{0}}$ is connected. Then $\cX$ is isomorphic to
  either a member of family 1, a member of family 2, or the
  premaniplex $\cX_s$.
\end{prop}

As discussed before, a semiregular facet-transitive 3-maniplex with
trivial facet stabilizer and regular symmetry type graph must have a
premaniplex $\cX$ satisfying the conditions of \cref{prop:regpmpx} as
its symmetry type graph.  However, we will see that it is impossible
for a 3-maniplex $\cM$ with trivial facet stabilizer to have $\cX_s$
as its symmetry type graph.

The \emph{$k$-gonal hosohedron} ($k \in \N \cup \{\infty\}$, $k\geq 2$),
denoted by $\{2,k\}$ is the regular polyhedron consisting of two
vertices joined by $k$ edges forming $k$ digons. The \emph{$2k$-gonal
  hemi-hosohedron} ($k\geq 1$ finite), denoted by $\{2,2k\}/2$, is the
quotient of the $2k$-gonal hosohedron by identifying opposite
flags. One might think of $\{2,2k\}/2$ as a map on the projective
plane with a single vertex connected to itself by $k$ edges that
touch it on opposite sides. In particular $\cX_s = \{2,2\}/2$.

\begin{lema}\label{lema:DigonExt}
  Let $\cM$ be a 3-maniplex with only digonal faces. Then $\cM$ is
  either a hosohedron or a hemi-hosohedron. If $\cM$ is a
  polyhedron, it must be a hosohedron.
\end{lema}
\begin{dem}
  Since facets are digons, the monodromy $r_0$ commutes with $r_1$,
  and since it commutes with $r_2$, it must be in the center of
  $\Mon(\cX)$. This implies that there are at most two vertices, which
  are the orbits under $\gen{r_1,r_2}$ of $\Phi$ and $\Phi^0$. If
  these two orbits are different (which is necessary for $\cM$ to be a
  polyhedron), then $r_0$ maps any flag in one orbit to a flag in the
  other one and $\cM$ is a hosohedron. If there is only one vertex, we
  may rename the flags as $\{1,2,\ldots,2k\}$ in such a way that $r_1$
  is the permutation $(1,2)(3,4)\cdots(2k-1,2k)$ and $r_2$ is the
  permutation $(2,3)(4,5)\cdots(2k,1)$. Notice that for every
  $i\in\{1,2,\ldots,2k\}$, one element of $\{r_1,r_2\}$ maps $i$ to
  $i+1$, while the other one maps $i$ to $i-1$. Since $r_0$ commutes
  with both $r_1$ and $r_2$ we get that $r_0(i+1)= r_0(i)\pm 1$. It
  follows by induction that either $r_0(i)=i+j$ for some fixed $j$ and
  for all $i$ or $r_0(i)=j-i$ for some fixed $j$ and for all
  $i$. Suppose $r_0(i)=j-i$ for all $i$. If $j$ is even, then
  $r_0(j/2)=j/2$, a contradiction since $r_0$ cannot have fixed points
  in a maniplex. If $j$ is odd, then
  $r_0((j-1)/2) = (j+1)/2 = (j-1)/2+1$, but this is either
  $r_1((j-1)/2)$ or $r_2((j-1)/2)$, a contradiction since $r_0r_1$ and
  $r_0r_2$ cannot have fixed points in a maniplex. Therefore,
  $r_0(i)=i+j$ for all $i$ and a fixed $j$. Finally, since $r_0$ must
  have order 2, $j=k$ and $\cM$ is a hemi-hosohedron.

  \qed
\end{dem}

Since hosohedra and hemi-hosohedra are all regular,
\cref{lema:DigonExt} implies that all maniplexes with only digonal
faces are regular. In particular, there are no maniplexes with trivial
facet stabilizer with symmetry type graph $\cX_s$. Therefore
\cref{prop:regpmpx} implies the following proposition:

\begin{prop}\label{prop:regsymtype}
  Let $\cM$ be a semiregular facet-transitive 3 maniplex with trivial
  facet stabilizer. Let $\cX$ be the symmetry type graph of
  $\cM$. Then, if $\cX$ is regular it is isomorphic to a premaniplex in
  either family 1 or family 2.
\end{prop}

\section{Constructing polyhedra with symmetry type graph in family 1}\label{sec:fam1const}

\begin{figure}[h!]
  \centering
  \includegraphics[width=0.5\linewidth]{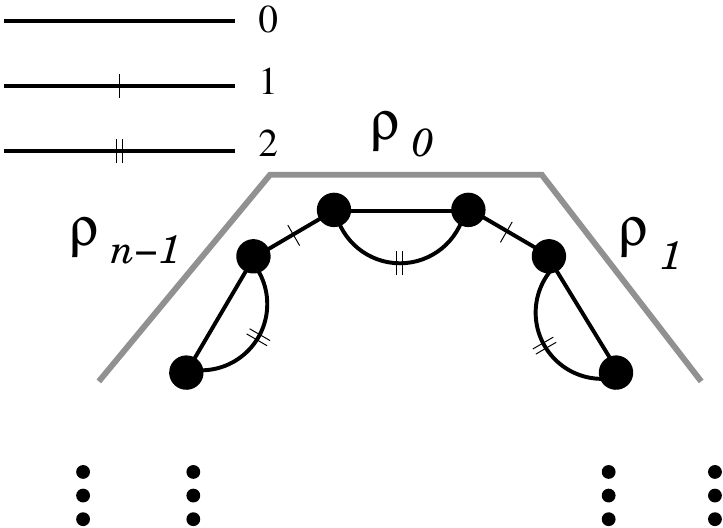}
  \caption{A premaniplex of family 1 with a voltage assignment.}
  \label{fig:fam1volts}
\end{figure}

Now we want to assign voltages to the premaniplexes of the previous
section to get examples of maniplexes and polytopes with trivial facet
stabilizer and with these symmetry type graphs. First, let $\cX$ be
the premaniplex of family 1 with $n$ sides. We want to find a voltage
assignment $\xi$ on $\cX$ such that $\cX^\xi$ has symmetry type $\cX$
and $n$-gonal faces. To ensure that $\cX^\xi$ has $n$-gonal faces it
is enough to assign trivial voltage to all the edges of color
different than 2. Now let $\alpha_0, \alpha_1,\ldots, \alpha_{n-1}$ be
the voltages of the edges of color 2 in cyclical order. To make sure
that $\cX^\xi$ is a maniplex we only need to make sure that $\alpha_i$
has order exactly 2 for every $i$. And finally, to make sure that
$\cX^\xi$ has symmetry type $\cX$, we only need to make sure that no
automorphism of $\cX$ induces an automorphism of our voltage group
$\Gamma=\gen{\alpha_0,\ldots,\alpha_{n-1}}$. Given a base flag $x$ in
$\cX$, we can can find a closed path $W_i$ whose voltage is $\alpha_i$
by using edges of colors 0 and 1 to get to the starting point of the
dart with voltage $\alpha_i$, then using that dart and going back to
the base using colors 0 and 1.  Given an automorphism $\tau$ of $\cX$,
the voltage of $W_i\tau$ will be the voltage of the image of the
dart with voltage $\alpha_i$. So, if $\tau$ induces an automorphism of
$\Gamma$, this automorphism must be an extension of the permutation
that $\tau$ induces on the generators
$\{\alpha_0,\ldots,\alpha_{n-1}\}$.  One way to make sure that no
non-trivial automorphism of $\cX$ induces an automorphism of $\Gamma$
is by using the automorphism group of a regular $n$-polytope $\cP$ as
$\Gamma$ and making $\alpha_i = \rho_i$, where $\rho_i$ is the
distinguished generator of $\Gamma(\cP)$ that maps a base flag $\Phi$
of $\cP$ to its $i$-adjacent flag $\Phi^i$. If $\cP$ is
non-degenerate, then the only possible non-trivial permutation of
distinguished generators that can induce an automorphism of
$\Gamma(\cP)$ is the one mapping $\rho_i$ to $\rho_{n-1-i}$ for every
$i$, as any other permutation would map a pair of non-commuting
generators to a pair that commutes. This permutation corresponds to
some reflection of $\cX$ and it induces an automorphism of the group
if and only if $\cP$ is self-dual. So, if $\cP$ is not self-dual and
non-degenerate, $\cX^\xi$ is a maniplex with symmetry type $\cX$.
This voltage assignment has been illustrated in \cref{fig:fam1volts},
where the label of an edge of the polygon corresponds to the voltage
of the darts of color 2 that start on that edge.

Finally, we would like for $\cX^\xi$ to be a polyhedron
(3-polytope). We will use the following theorem:

\begin{teo}\label{teo:PolyPreExt}\cite[Proposition 3.14]{CayExt}
  Let $(\cX,\Gamma,\xi)$ be a voltage $(r+1)$-premaniplex. Suppose
  $\cX_{\ol{r}}$ is connected, $\xi(d)=1$ for every dart $d$ of color
  smaller than $r$ and that $\cX^\xi$ is a maniplex. Then $\cX^\xi$ is
  a polytope if and only if $\cX_{\ol{r}}$ is a polytope, and whenever
  there is a path $W$ with voltage 1 and colors in $[m,r]$, there is
  also a path $W'$ with the same endpoints as $W$ that uses only
  colors in $[m,r-1]$.
\end{teo}

In our case $r=2$. Notice that by hypothesis $\cX_{\ol{r}}$ is
connected, which makes the case $m=0$ trivial. On the other hand, for
$m=2$ one only has to notice that there are no edges of color 2 with
trivial voltage. Therefore, the only interesting case is when
$m=1$. In other words, we want to make sure that if two flags are
connected by a path with colors 1 and 2 and voltage 1, then they are
either the same flag or 1-adjacent flags.

Let us start by analyzing the voltage of an alternating path $W$ of
colors 1 and 2 starting from $\Phi$. Suppose $W$ uses a total of $k$
darts of color 2. We can see that if $W$ starts with a dart of color
2, then
\[
 \xi(W) = \prod_{i=0}^{k-1}\rho_i
\]
where the product goes from right to left and the subscripts are taken
modulo $n$. On the other hand, if the first dart in $W$ has color 1,
then
\[
 \xi(W) = \prod_{i=0}^{k-1}\rho_{n-1-i}.
\]
However, if $\prod_{i=0}^{n-1}\rho_i$ has order $\ell$, we might write
this as
\[
 \prod_{i=0}^{k-1}\rho_{n-1-i}=\prod_{i=0}^{n\ell-k-1}\rho_i.
\]
Notice that if the left-hand side has $k$ factors, the right-hand side
has $n\ell-k$. Regardless of the color of the first dart, the path $W$
ends in $\Phi$ or $\Phi^1$ if and only if $k$ is a multiple of
$n$. Therefore, what we want 
is that
$ \prod_{i=0}^{k-1}\rho_i = 1$ if and only if $k$ is a multiple of
$n$, and the same is true for
$\xi(W) = \prod_{i=0}^{k-1}\rho_{n-1-i}$.

This condition can be interpreted in terms of the Petrie polygon of
$\cP$. Typically, a Petrie polygon of a polytope is defined as a
cyclical sequence of edges (1-faces) such that whenever
$2\leq i \leq n-1$, any $i$ consecutive edges share an $i$-face but
any $i+1$ consecutive edges do not. According to this definition some
degenerate polytopes and some non-polytopal maniplexes would not have
Petrie polygons, so the following alternative is preferred: We will
define a \emph{Petrie polygon} as a flag orbit under the action of the
monodromy $\omega_\pi:=r_{n-1}r_{n-2}\cdots r_0$. 
The cardinality $\ell$ of a Petrie polygon $P$ is called its
\emph{length}. If we look at the sequence of edges
$\{(\omega_\pi^i\Phi)_1\}_{i=0}^\ell$ on a non-degenerate polytope, we
would see that it satisfies the usual definition of Petrie
polygon. Actually, in~\cite{PetrieSchemes} the sequence
$\{(\omega_\pi^i\Phi)_1\}_{i=0}^\ell$ is called a \emph{Petrie scheme}
and its sequence of edges is called a Petrie polygon, but in this
article we do not make that distinction. Notice that the edge
$(\omega_\pi^i\Phi)_1$ is incident to the vertices
$(\omega_\pi^i\Phi)_0$ and $(\omega_\pi^{i+1}\Phi)_0$. In this
definition we make sure that every flag is in exactly one Petrie
polygon. 
If any two different flags in $P$ have different vertices (0-faces),
we will say that $P$ is \emph{simple}. 

The following lemma will be useful:

\begin{lema}\label{lema:PetrieSimple}
  Let $\cP$ be an $n$-polytope. Let $\Phi$ be a flag of $\cP$ and let
  $P$ be the Petrie polygon containing $\Phi$. Let $\ell$ be the
  length of $P$ and let $S=\{\Psi_i\}_{i=0}^{n\ell-1}$ be the sequence
  of flags defined by $\Psi_0=\Phi$ and $\Psi_{i+1} = r_i \Psi_i$ for
  $0\leq i \leq n\ell$, where the subscript of $r_i$ is taken modulo
  $n$ and the subscript of $\Psi_i$ modulo $n\ell$. If $P$ is simple,
  then all the elements of $S$ are different.
\end{lema}

\begin{dem}
  Note that $P=\{\Psi_{ni}\}_{i=0}^{\ell-1}$. Notice that whenever
  $\floor{(j-1)/n}=i-1$ we get that $\Psi_j$ has the same vertex
  (0-face) as $\Psi_{ni}$. Therefore, since $P$ is simple, we get that
  $\Psi_i$ and $\Psi_j$ have the same vertex if and only if
  $\floor{(i-1)/n} = \floor{(j-1)/n}$. As a small caveat, for this
  statement to be completely true we must rename $\Psi_0$ as
  $\Psi_{n\ell}$.

  Let $i$ and $j$ be such that $\Psi_i=\Psi_j$ with
  $0\leq i < j\leq n\ell-1$.
  Then
  $\Psi_{j}=\prod_{k=i}^{j-1}r_k \Psi_{i}$. Since $\Psi_{i}$ and
  $\Psi_{j}$ have the same vertex, we get that
  $\floor{(i-1)/n} = \floor{(j-1)/n}$. 
  In particular $j-i< n$, so $r_{i}$ appears exactly once in the
  product $\prod_{k=i}^{j-1}r_k$. But since $\cP$ is a polytope, this
  implies that $\prod_{k=i}^{j-1}r_k$ must change the $i$-face of
  every flag, so we have reached a contradiction.
  \qed
\end{dem}

Going back to the case where $\cP$ is a regular $n$-polytope with
simple Petrie polygons, using that $\Phi\rho_i = r_i \Phi$ and the
fact that automorphisms commute with monodromies one can see that
\[
  \Phi\prod_{i=0}^k\rho_i = \prod_{i=0}^kr_i \Phi = r_k
  \prod_{i=0}^{k-1} r_i\Phi.
\]
So the sequence $S$ in \cref{lema:PetrieSimple} is the same as
$\{\Phi\prod_{i=0}^{k-1}\rho_i\}_{k=0}^{n\ell-1}$ where the subscripts are
taken modulo $n$. Applying the lemma we get that
$\Phi\prod_{i=0}^{k-1}\rho_i = \Phi$ if and only if $k$ is divisible
by $n\ell$. In particular, if $\prod_{i=0}^{k-1}\rho_i = 1$ then $k$
is divisible by $n$.

Now let $W$ be a path with darts of colors 1 and 2. Without loss of
generality, we may assume that $W$ is reduced, that is, that it does
not have any pair of consecutive inverse darts. This is equivalent to
saying that $W$ alternates between colors 1 and 2. Let $W'$ be the
path that follows the same color sequence as $W$ but starting at the
base flag of $\cX$. Notice that the voltage of $W'$ is a conjugate of
either the voltage of $W$ or the voltage of $W^{-1}$ (depending on
whether $W$ starts at a left or a right flag), so $\xi(W)=1$ if and
only if $\xi(W')=1$. But we have proved that $\xi(W')=1$ if and only
if $W'$ ends at either the base flag or its 1-adjacent flag, and since
$\cX$ is regular, the same must be true for $W$.

Therefore, we have arrived to the following result:

\begin{teo}\label{teo:fam1const}
  Let $\cP$ be a regular $n$-polytope with simple Petrie
  polygons. Suppose that $\cP$ is non-degenerate and
  non-self-dual. Let $(\cX,\Gamma,\xi)$ be the voltage premaniplex
  where $\cX$ is the premaniplex of family 1 with $n$ sides,
  $\Gamma$ is the automorphism group of $\cP$ with distinguished
  generators $\{\rho_0,\ldots,\rho_{n-1}\}$, and $\xi$ is the voltage
  assignment that assigns voltage $\rho_i$ to the darts starting at
  the $i$-th 1-face (as in \cref{fig:fam1volts}). Then $\cX^\xi$ is a
  polyhedron with symmetry type graph isomorphic to $\cX$. In
  particular $\cX^\xi$ is vertex-transitive, facet-transitive and it
  has trivial facet stabilizer.
\end{teo}

\begin{ej}
  Let $\cP$ be a 3-dimensional cube, and let $\rho_0$, $\rho_1$ and
  $\rho_2$ be the distinguished generators for its automorphism
  group. Let $\cX$ be the premaniplex in family 1 with 3 sides. Then,
  the derived graph from the voltage graph in \cref{fig:fam1volts} is
  a semiregular facet-transitive polyhedron with 48 triangular
  faces, 72 edges and 8 vertices of degree 18. It has genus 9 and it
  has trivial facet stabilizer. If we replace the cube
  with the octahedron we get the same resulting polyhedron.

  Analogously, both the icosahedron and the dodecahedron give us a
  polyhedron with 120 triangular faces, 180 edges and 12 vertices of
  degree 30. It has genus 25 and trivial facet stabilizer.

  Finally, if $\cP$ is the tetrahedron the result is a polyhedron with
  24 triangular faces, 36 edges and 6 vertices of degree 12. It has
  genus 4 and its facet stabilizer has order two, as the tetrahedron
  is self-dual.
\end{ej}

Notice that the degree of the vertices in $\cO(\cP)$ is $n$ times the
length of the Petrie polygons of $\cP$.

\section{Constructing polytopes with symmetry type graph in family 2}\label{sec:fam2const}

To construct polytopes with symmetry type graphs in family 2 we
proceed in a similar fashion. Let $\cX$ be the premaniplex in family 2
constructed from a $4k$-gon (see \cref{fig:fam2}). Once again, we
enumerate the 1-faces of $\cX$ so that $e_i$ is the 1-face of the flag
$(r_1 r_0)^i \Phi$ where $\Phi$ is the base flag. Note that an edge of
color 2 joins a flag in the 1-face $e_i$ with one in the edge
$e_{2k+i}$ (with the subscripts taken modulo $4k$). We assign trivial
voltage to the darts of colors 0 and 1 to ensure that the facets of
the derived graph are $4k$-gons. If we assign a voltage $\alpha_i$ to
a dart going from edge $e_i$ to edge $e_{2k+i}$, we must assign the
same voltage to the other dart that goes from $e_i$ to
$e_{2k+i}$. This is to ensure that the derived graph is a maniplex.

\begin{figure}[h!]
  \centering
  \includegraphics[width=0.5\linewidth]{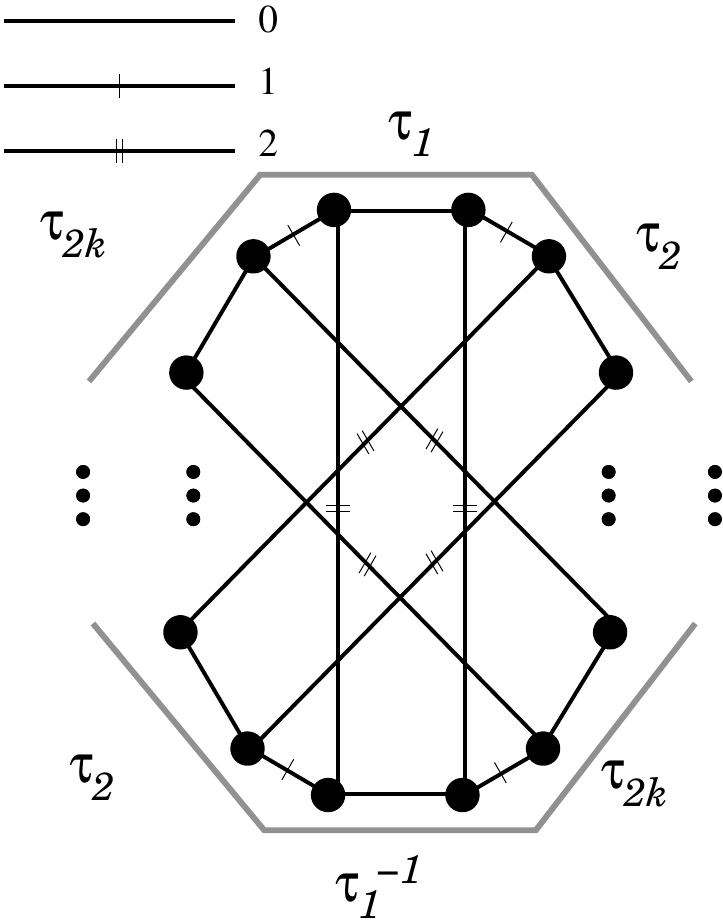}
  \caption{A premaniplex of family 2 with a voltage assignment.}
  \label{fig:fam2volts}
\end{figure}

To assign the voltages to the darts of color 2, this time we use the
automorphism group of a chiral polytope. A maniplex is \emph{chiral}
if it has 2 orbits on flags and adjacent flags are in different
orbits. In 2010 Pellicer~\cite{Quirales} proved that there are
abstract chiral polytopes in every rank greater or equal to 3. Let
$\cP$ be a chiral polytope of rank $2k+1$ and let $\Phi$ be a base
flag of $\cP$. Because of chirality, for every $i$ in
$\{1,2,\ldots,2k\}$, there is an automorphism $\tau_i$ that maps
$\Phi$ to $r_ir_0\Phi$. It is worth noting that $\tau_i$ is an
involution if and only if $i\neq 1$ and $\tau_i\tau_j$ is an
involution if and only if $|i-j|>1$. This is a consequence of the
known fact that there are no degenerate chiral
polytopes~\cite[Lemma 2.3.16]{ACP}. For $i\in\{0,1,\ldots,2k-1\}$ we
will assign the voltage $\tau_{i+1}$ to the darts of color 2 in $\cX$
that connect the 1-face $e_i$ to $e_{2k+i}$. This voltage assignment
is illustrated in \cref{fig:fam2volts}, where once again, the labels
of the edges of the polygon correspond to the voltages of the darts of
color 2 that start on that edge.

Using the fact that $\tau_1$ is the only $\tau_i$ that is not an
involution we can see that the only possible automorphisms of $\cX$
that could induce an automorphism of $\Gamma(\cP)$ are the identity,
the half turn, the reflection that fixes $e_0$ and
$e_{2k}$, and the reflection that fixes $e_k$ and $e_{3k}$.

It is known (see~\cite[Corollary~2.5]{ExtQuiralesToros}) that if $\cP$
is a chiral polytope, there is no automorphism of $\Gamma(\cP)$ that
maps $\tau_1$ to $\tau_1^{-1}$ and fixes every other
$\tau_i$. Therefore, using \cref{teo:lifts} we get that the half-turn
does not lift.

The reflection that fixes $e_0$ fixes edges with voltage $\tau_1$ and
for $2\leq i \leq 2k$, it swaps edges with voltage $\tau_i$ with edges
with voltage $\tau_{2k+2-i}$. Notice that $\tau_1\tau_2$ is not an
involution but $\tau_1\tau_{2k}$ is whenever $k\geq 2$. Therefore, the
reflection that fixes $e_0$ does not lift for $k\geq 2$. Notice that
for $k=1$, the reflection that fixes $e_0$ induces the identity
automorphism in $\Gamma(\cP)$, so in that case it does lift and the
derived maniplex has non-trivial facet stabilizer (the stabilizer
would be the group of order of order 2 generated by the lift of this
reflection).

The reflection that fixes $e_k$ swaps edges with voltage $\tau_1$ with
edges with voltage $\tau_1^{-1}$, and it swaps edges with voltage
$\tau_i$ with edges with voltage $\tau_{2k+2-i}$. For $k=1$, this
reflection cannot induce an automorphism of $\Gamma(\cP)$, since it
would fix $\tau_2$ while inverting $\tau_1$. For $k\geq 2$, we see
that $\tau_1\tau_2$ is not an involution while $\tau_1^{-1}\tau_{2k}$
is, so once again this reflection does not lift.

We have proved that for $k\geq 2$ no non-trivial automorphism of $\cX$
lifts, and therefore $\cX$ is the symmetry type graph of $\cX^\xi$.


Now we want to see when our examples are polytopes.

If we look at a path $W$ of alternating colors 2 and 1 (starting with
color 2) that starts on the left flag of $e_0$ and uses $j$ edges of
color 2, we will see that its voltage is
$\prod_{i=1}^j \tau_i^{\epsilon_i}$, where the subscript of $\tau_i$
is taken modulo $2k$ (using $\tau_{2k}$ when $i\equiv 0 \mod 2k$) and

\begin{equation}
  \label{eq:epsilon}
  \epsilon_i=
  \begin{cases*}
    1 &  if  $i\not\equiv 2k+1\mod 4k$\\
    -1 & if $i\equiv 2k+1\mod 4k.$
  \end{cases*}
\end{equation}
The path $W$ ends in
$\Phi$ or $\Phi^1$ whenever it uses a multiple of $4k$ edges of color
2. So using \cref{teo:PolyPreExt}, we get that in order to find
examples where the derived graph is polytopal, we need to make sure
that $\prod_{i=1}^j \tau_i^{\epsilon_i}=1$ implies that $j$ is
divisible by $4k$.

\begin{prop}\label{prop:prod_tau}
  Let $\cP$ be a chiral $(2k+1)$-polytope with base flag $\Phi$, and
  let $\tau_i$ be such that $\Phi\tau_i=r_ir_0 \Phi$ for
  $1\leq i \leq 2k$ and let $\epsilon_i$ be defined as in
  \cref{eq:epsilon}. Let $j$ be a positive integer. If
  $j\not\equiv \floor{(j-1)/(2k)} \mod 2$, then
  \[
    \Phi \prod_{i=1}^j \tau_i^{\epsilon_i} =
    \prod_{i=0}^{j+\floor{(j-1)/(2k)}} r_i \Phi,
  \]
  where the subscript of $r_i$ is taken modulo $2k+1$ and the
  subscript of $\tau_i$ modulo $2k$.  On the other hand, if
  $j\equiv \floor{(j-1)/(2k)} \mod 2$, then
  \[
    \Phi \prod_{i=1}^j \tau_i^{\epsilon_i} =
    r_0\prod_{i=0}^{j+\floor{(j-1)/(2k)}} r_i \Phi.
  \]
\end{prop}
\begin{dem}
  We will prove this in a ``hand-wavy'' manner. The reader may
  formalize this proof using induction and making several cases.
  For $j=1$ we get $\Phi\tau_1 = r_1 r_0 \Phi$ by definition of
  $\tau_1$. Now suppose the result is true for $j-1$.  Using the
  definition of $\tau_j$ and the fact that automorphisms commute with
  monodromies we get

  \begin{eqnarray*}
    \Phi \prod_{i=1}^j \tau_i^{\epsilon_i}
    &=&
        \Phi\tau_j^{\epsilon_j}\prod_{i=1}^{j-1} \tau_i^{\epsilon_i}\\
    &=&
        ((r_jr_0)^{\epsilon_j}\Phi)\prod_{i=1}^{j-1} \tau_i^{\epsilon_i}\\
    &=&
        (r_jr_0)^{\epsilon_j}(\Phi\prod_{i=1}^{j-1} \tau_i^{\epsilon_i}).\\
  \end{eqnarray*}
  Recall that the subscripts of $\tau_i$ and $\tau_j$ are taken modulo
  $2k$, and therefore, the subscript $j$ of $r_j$ is also taken modulo
  $2k$. If $j\not\equiv 1\mod 2k$, we see that the exponent
  $\epsilon_j$ is irrelevant and the $r_0$ gets canceled or not
  depending on whether $ \Phi\prod_{i=1}^{j-1} \tau_i^{\epsilon_i}$ has
  the factor $r_0$ on the left or not. Notice that for $j<2k+1$ there is
  alternation between these two cases, where the even terms are the
  ones with an $r_0$ on the left. In particular
  $\Phi\prod_{i=1}^{2k} \tau_i^{\epsilon_i} = r_0\prod_{i=0}^{2k}
  r_i\Phi$. Then for $j=2k+1$ we get
  \begin{eqnarray*}
    \Phi \prod_{i=1}^{2k+1} \tau_i^{\epsilon_i}
    &=& (r_1r_0)^{\epsilon_{2k+1}}(\Phi\prod_{i=1}^{2k}
        \tau_i^{\epsilon_i})\\
    &=& r_0r_1( r_0\prod_{i=0}^{2k}
        r_i\Phi) \\
    &=& r_0 \prod_{i=0}^{2k+2}r_i\Phi.
  \end{eqnarray*}
  Notice how at this point the alternation pauses momentarily and we
  get two consecutive products with an $r_0$ on the left. Notice also
  that the superscript of the product of $r_i$ increases by 2 rather
  than 1. After this point alternation resumes, but now the odd terms
  are the ones with an $r_0$ on the left. In particular
  $\Phi\prod_{i=1}^{4k} \tau_i^{\epsilon_i} = \prod_{i=0}^{4k+1}
  r_i\Phi$. Now for $j=4k+1$ we get:
  \begin{eqnarray*}
    \Phi \prod_{i=1}^{4k+1} \tau_i^{\epsilon_i}
    &=& (r_1r_0)^{\epsilon_{4k+1}}(\Phi\prod_{i=1}^{4k}
        \tau_i^{\epsilon_i})\\
    &=& r_1r_0(\prod_{i=0}^{4k+1}
        r_i\Phi) \\
    &=& \prod_{i=0}^{4k+3}r_i\Phi.
  \end{eqnarray*}
  This time we get two consecutive terms without an $r_0$ on the left
  and once again the superscript has increased by 2. This patterns
  repeats periodically: usually we alternate between terms with and
  without $r_0$ on the left and we increase the superscript by 1, but
  when we get to a term congruent with 1 mod $2k$ the alternation stops
  and the superscript increases by 2.\qed
\end{dem}

Let $\cP$ be a chiral polytope of rank $2k+1$ with base flag $\Phi$
with simple Petrie polygons. Similarly to the case for family 1, we
take a path $W$ in $\cX$ starting from the base flag, alternating
darts of color 2 and 1, starting with color 1 and using a total of $j$
darts of color 2. Using \cref{prop:prod_tau} we see that the voltage
$W$ acts as a $2t$-step rotation in the Petrie polygon $P$ whenever
$W$ uses $4kt$ edges of color 2.

Suppose that $\xi(W)=1$, in particular $\Phi=\Phi\xi(W)$. We want to
prove that $j$, the number of edges of color 2 used by $W$, is a
multiple of $4k$.

If $j\not\equiv\floor{(j-1)/(2k)} \mod 2$, then
$\Phi\xi(W) = \prod_{i=0}^{j+\floor{(j-1)/(2k)}}r_i \Phi$. Let
$t=j+\floor{(j-1)/(2k)}+1$. We can use the same argument as we used
for family 1 to prove that $\Phi=\prod_{i=0}^{t-1} r_i \Phi$ implies
that $t$ is divisible by the rank, which in this case is $2k+1$.  A
simple substitution shows that when $j=2km$ for some positive $m$,
then $t=(2k+1)m$. Since $t$ is defined by a monotone function on
$j$, the converse is also true. This shows that if
$\prod_{i=0}^{t-1} r_i$ fixes the base flag, then $j$ is divisible by
$2k$. Since $j\not\equiv \floor{(j-1)/2k} \mod 2$, then
$j/(2k) \equiv j\equiv 0 \mod 2$, and therefore $j$ must be a multiple
of $4k$.

Now suppose $j\equiv \floor{(j-1)/2} \mod 2$ and
$\Phi=\Phi\xi(W)=r_0\prod_{i=0}^{j+\floor{(j-1)/(2k)}} r_i \Phi$. Then
$r_0\Phi=\prod_{i=0}^{j+\floor{(j-1)/(2k)}}r_i \Phi$. Let
$t=j+\floor{(j-1)/(2k)}+1$ and let $t'$ be its residue modulo
$(2k+1)\ell$, where $\ell$ is the length of $P$. Then, by definition
of $\ell$ we see that
$r_0\Phi=\prod_{i=0}^{t-1}r_i \Phi
=\prod_{i=0}^{t'-1}r_i\Phi$. \cref{lema:PetrieSimple} tells us that
$t'=1$, in particular $t\equiv 1 \mod 2k+1$. We established earlier that
$t$ is divisible by $2k+1$ if and only if $j$ is divisible by $2k$. If
$j\equiv r \mod 2k$, for $1\leq r<2k$, a substitution shows that
$t\equiv r+1 \mod 2k+1$. Therefore, it is not possible that
$t\equiv 1 \mod 2k+1$. This proves that this case cannot happen.

We have shown that $\xi(W)=1$ implies that $W$ ends at the base flag,
or at its 1-adjacent flag. In analogy to the case for family 1, we can
use this to prove that any path with colors 1 and 2 with trivial
voltage must end either on the flag where it started or on its
1-adjacent flag. This proves the following theorem:

\begin{teo}\label{teo:fam2const}
  Let $k\geq 2$ and let $\cP$ be a chiral $(2k+1)$-polytope with
  simple Petrie polygons. Let $(\cX,\Gamma,\xi)$ be the voltage
  premaniplex where $\cX$ is the premaniplex of family 2 with $4k$
  sides, $\Gamma$ is the automorphism group of $\cP$ with
  distinguished generators $\{\tau_1,\ldots,\tau_{2k}\}$, and $\xi$ is
  the voltage assignment that assigns voltage $\tau_{i+1}$ to the
  darts starting at the $i$-th 1-face for $i<2k$ (as in
  \cref{fig:fam2volts}). Then $\cX^\xi$ is a polyhedron with symmetry
  type graph isomorphic to $\cX$. In particular $\cX^\xi$ is
  vertex-transitive, facet-transitive and it has trivial facet
  stabilizer.
\end{teo}

\begin{ej}\label{ej:Fam2}
  Let $G$ be the group with presentation
  \begin{align*}
    \gen{\sigma_1,\sigma_2,\sigma_3,\sigma_4&| \mathcal{R}(5),
                                              \sigma_1^{12} = \sigma_2^6= \sigma_3^6 = \sigma_4^3 =
                                              [\sigma_1,\sigma_2^2] = [\sigma_4,\sigma_3^2] =
                                              (\sigma_2^{-1}\sigma_3)^4 \\
                                            &=
                                              \sigma_1^{-1}\sigma_2\sigma_3^3\sigma_1^{-1}\sigma_3^2\sigma_2\sigma_3^{-1}= 1}
  \end{align*}
  where $\mathcal{R}(5)$ is a set of relations satisfied by every
  chiral 5-maniplex (see~\cite[Section 4]{AbelianCoversChirals}). In
  \cite[Example 4.20]{AbelianCoversChirals} it is shown that $G$ is
  the automorphism group of a chiral 5-polytope $\cP$ with 41472
  flags. The generator $\sigma_i$ maps the base flag $\Phi$ to
  $\Phi^{i,i-1}$, and one can show that if $\tau_i$ is defined as the
  automorphism of $\cP$ that maps $\Phi$ to $\Phi^{0,i}$ then
  $\sigma_i = \tau_{i-1}\tau_i^{-1}$ ($\tau_0$ is the identity) and
  $\tau_i = \sigma_i^{-1}\sigma_{i-1}^{-1}\cdots \sigma_1^{-1}$. Using
  the package RAMP~\cite{RAMP0.7.01} for the software
  GAP~\cite{GAP4.12.1}, we can verify that $\cP$ has simple Petrie
  polygons of length 12.

  If $\cX$ is the premaniplex in family 2 with 8 sides, the derived
  graph $\cX^\xi$ from the voltage graph in \cref{fig:fam2volts} is a
  5-polytope with 20736 octagonal faces, 82944 edges and 3456
  vertices of degree 48. It has genus 29377.
\end{ej}

To calculate the degree of the vertices we have to remember that a
$4k$-step rotation (where $2k+1$ is the rank of $\cP$) around a vertex
of $\cX^\xi$ corresponds to a 2-step rotation on the Petrie polygon of
$\cP$.  So if $\ell$ is the length of the Petrie polygon containing
the base flag of $\cP$, the degree of the vertices of $\cX^\xi$ is
$2k\ell$ if $\ell$ is even, and $4k\ell$ if $\ell$ is odd. However, in
odd ranks the length of the Petrie polygons of chiral polytopes
cannotbe odd: this is because their existence would imply that there
is a closed path of odd length in the maniplex, but each closed path
must have the same number of flags in one orbit as in the other one.

If we were to prove that there are chiral $(2k+1)$-polytopes
($k\geq 2$) with simple Petrie polygons, \cref{teo:fam2const} would
imply the existence of polyhedra whose symmetry type graph is the
premaniplex of family 2 with $4k$ sides. So we ask the following
question

\begin{preg}
  Are there chiral polytopes with simple Petrie polygons in any (odd)
  rank $n\geq 3$?
\end{preg}

It might be enough to analyze the examples constructed by Pellicer
in~\cite{Quirales} or by Montero and Toledo
in~\cite{ExtQuiralesToros}.

\section{Interpretation via voltage operations}\label{sec:VoltOps}

In this section we want to interpret the constructions of the previous
two sections as applying some well defined operation to a regular or
chiral polytope. To do this we will use voltage operations. Voltage
operations were introduced in~\cite{VoltOps}, and they are a family of
operations that include the Wythoffian constructions, duality, Petrie
duality, and even some operations that modify the rank, such as the
prism and the pyramid. They are very useful, because they can be
applied not only to maniplexes and polytopes, but also to
premaniplexes. In fact, if $\cO$ is a voltage operation, $\cM$ is a
maniplex and $\cX$ is its symmetry type graph, then the symmetry type
graph of $\cO(\cM)$ is a quotient of $\cO(\cX)$ (see~\cite[Theorem
5.1]{VoltOps}).

An \emph{$(n,m)$-voltage operator} is a voltage premaniplex
$(\cY, W^n, \eta)$ where $\cY$ has rank $m$ and $W^n$ is the Coxeter
group
$\gen{r_0,r_1,\ldots r_{n-1}\mid r_i^2=(r_ir_j)^2=1 \text{ whenever }
  |i-j|>1}$. Since the monodromy group of any $n$-premaniplex is a
quotient of $W^n$, this group acts naturally on the set of flags of
any $n$-premaniplex (on the left) by monodromies. Note that we are
abusing notation by using $r_i$ both for an element of $W^n$ and an
element of $\Mon(\cX)$ where $\cX$ is a premaniplex, however this is
not a problem because both elements act the same way on each flag of
$\cX$.

Given an $n$-premaniplex $\cX$ and an $(n,m)$-voltage operator
$(\cY,W^n,\eta)$ we can create a new $m$-premaniplex denoted by $\cX
\ertimes \cY$ as follows: The flags of $\XY$ are the Cartesian product
$\cX\times\cY$, and the $i$-adjacent flag of $(x,y)$ is given by
\[
  r_i(x,y) = (\eta(^iy)x,r_i(y)),
\]
where $^iy$ is the dart of color $i$ that starts at $y$. In general,
$\XY$ might not be connected. However, if $\cY$ has a spanning tree
all of whose darts have trivial voltage, and the voltages of the
remaining darts generate $W^n$, then $\XY$ is connected
(see~\cite[Corollary3.10]{VoltOps}).

Given a voltage $n$-premaniplex $(\cX,\Gamma,\xi)$ and a voltage
operator $(\cY,W^n,\eta)$,~\cite[Theorem 5.4]{VoltOps} gives us a way
to construct a voltage assignment $\theta$ in $\XY$ with values in
$\Gamma$ such that $(\XY)^\theta$ is isomorphic to $\cX^\xi\ertimes
\cY$. Given a monodromy $w=r_{i_1}r_{i_2}\cdots r_{i_\ell}$ let us
define the path $^wx$ as the path that starts in $x$ and follows the
colors $i_\ell,\ldots,i_2,i_1$.

\begin{teo}\label{teo:theta}\cite[Theorem 5.4]{VoltOps}
  Let $(\cX,\Gamma,\xi)$ be a voltage $n$-premaniplex and
  $(\cY,W^n,\eta)$ an $(n,m)$-voltage operator. Define
  \[
    \theta(^i(x,y)):=\xi(^{\eta(^iy)}x).
  \]
  Then, $(\XY)^\theta$ is isomorphic to $\cX^\xi\ertimes \cY$.
\end{teo}

As a corollary we get
\begin{coro}\cite[Corollary 5.6]{VoltOps}\label{coro:OpsRegs}
  Let $\cM$ be a regular maniplex with automorphism group
  $\gen{\rho_0,\rho_1,\ldots,\rho_{n-1}}$ and let $(\cY,W^n,\eta)$ be
  a voltage operator. Then $\MY$ is isomorphic to $\cY^\nu$ where
  $\nu(d)$ is obtained by writing $\eta(d)$ as a word in
  $\{r_0,\ldots,r_{n-1}\}$ and replacing every occurrence of $r_i$
  with the corresponding $\rho_i$.
\end{coro}

\subsection{Family 1}\label{sec:VoltOps1}
\begin{figure}[h!]
  \centering
  \includegraphics[width=0.5\linewidth]{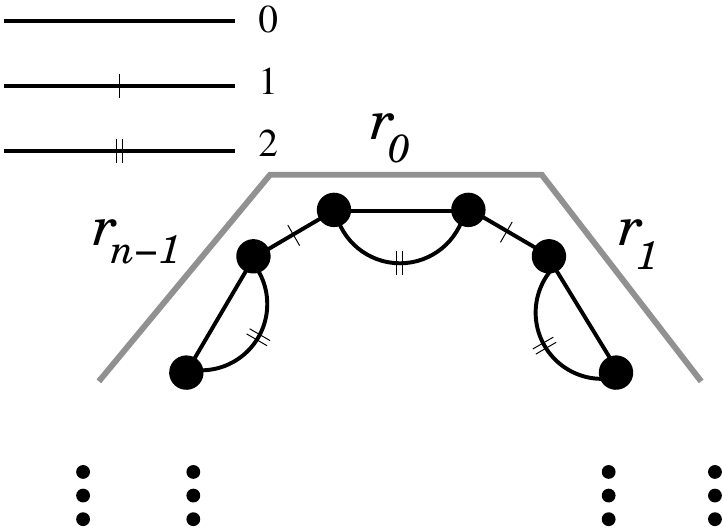}
  \caption{The voltage operator $(\cX,W^n,\eta)$.}
  \label{fig:fam1op}
\end{figure}

Let $\cX$ be the member of family 1 obtained from an
$n$-gon. \cref{teo:fam1const} tells us that we can construct a
polyhedron whose symmetry type is $\cX$ by using the voltage
premaniplex in \cref{fig:fam1volts}, where
$\gen{\rho_0,\ldots,\rho_{n-1}}$ is the automorphism group of some
non-degenerate, non-self-dual regular $n$-polytope $\cP$ with simple
Petrie polygons. Then, using \cref{coro:OpsRegs} we get that $\cX^\xi$
is isomorphic to $\cP\ertimes \cX$ where $(\cX,W^n,\eta)$ is the
voltage operation we get by replacing the voltage $\rho_i$ by $r_i$
(see \cref{fig:fam1op}).

\begin{figure}[h!]
  \centering
  \includegraphics[width=0.7\linewidth]{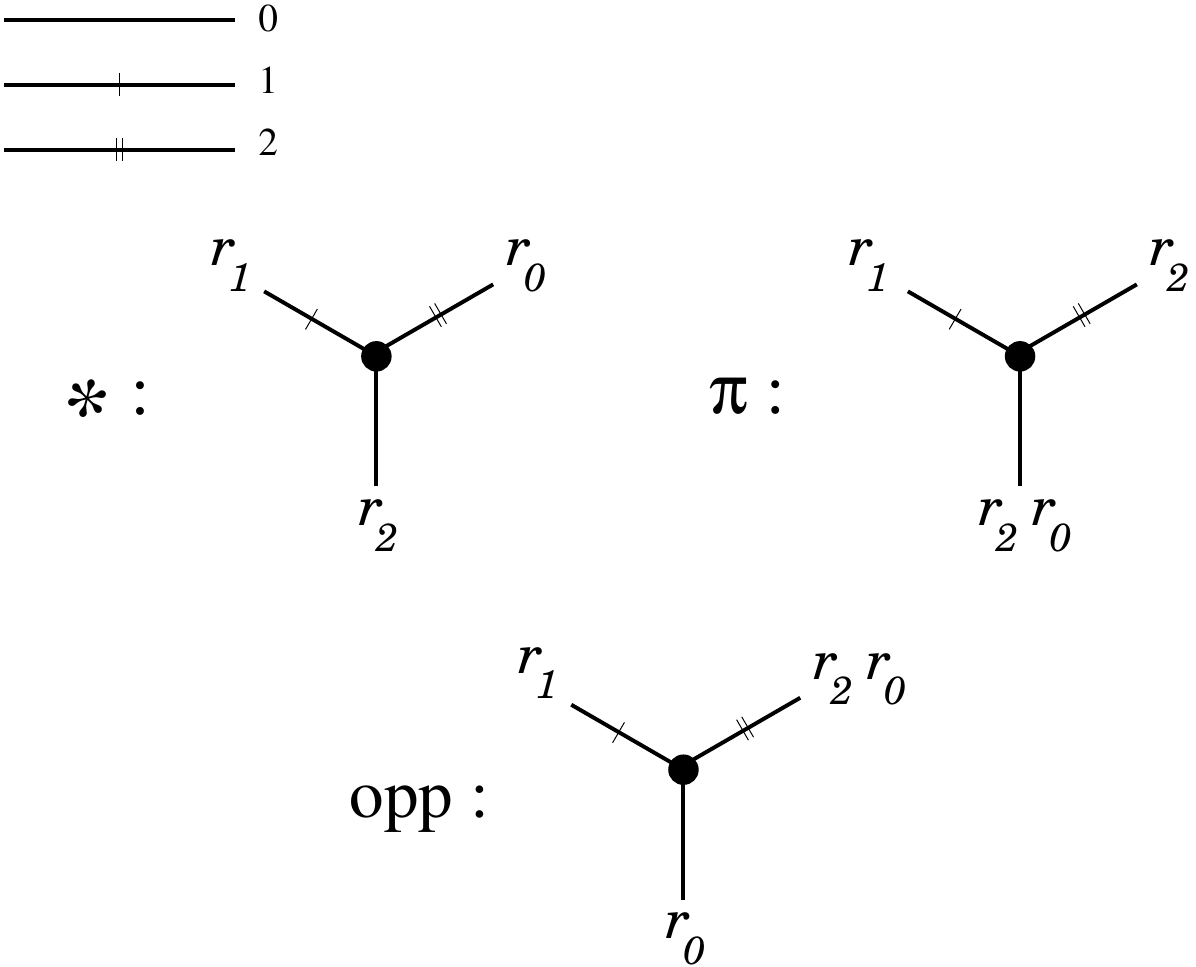}
  \caption{The voltage operators for the dual, Petrial, and opposite
    operations.}
  \label{fig:1vert}
\end{figure}

In what follows we want to give some topological interpretation of
what the voltage operation $\cO:\cP\mapsto \cP\ertimes \cX$
does. First, we need to introduce more concepts. Given a polyhedron
$\cP$ with simple Petrie-polygons, we define the \emph{Petrial} or
\emph{Petrie dual} of $\cP$, as the polyhedron $\cP^\pi$ with the same
vertices and edges, but where the 2-faces correspond to the Petrie
polygons of $\cP$. We can extend this definition to all 3-maniplexes
(and pre-maniplexes) as follows: Given a 3-premaniplex $\cX$, its
Petrial $\cX^\pi$ is defined as the premaniplex we get by removing the
edges of color 0, and adding a new edge of color 0 between each flag
$\Phi$ and its $02$-adjacent flag $r_2r_0\Phi$. If we denote by
$\cX^*$ the usual dual of $\cX$, then $((\cX^*)^\pi)^*$ is denoted by
$\opp{\cX}$ and is known as the \emph{opposite} of $\cX$. It is easy
to see that the operations dual, Petrial and opposite are voltage
operations. The corresponding voltage operators are illustrated in
\cref{fig:1vert}.

Now we want to use the version of \cref{teo:theta} for the
concatenation of two voltage operations:

\begin{teo}\label{teo:thetacomp}\cite[Theorem 6.1]{VoltOps}
  Let $\cX$ be an $n$-premaniplex, $(\cY_1,W^n,\eta_1)$ an
  $(n,m)$-voltage operator and $(\cY_2,W^m,\eta_2)$ an
  $(m,\ell)$-voltage operator. Then
    \[
      (\cX\rtimes_{\eta_1}\cY_1)\rtimes_{\eta_2}\cY_2 \cong
      \cX\rtimes_\theta(\cY_1\rtimes_{\eta_2}\cY_2),
    \]
    where
    \[
      \theta(^i(y_1,y_1)) = \theta(^{\eta_2(^i y_2)} y_1).
    \]
  \end{teo}

\begin{figure}[h!]
  \centering
  \includegraphics[width=0.5\linewidth]{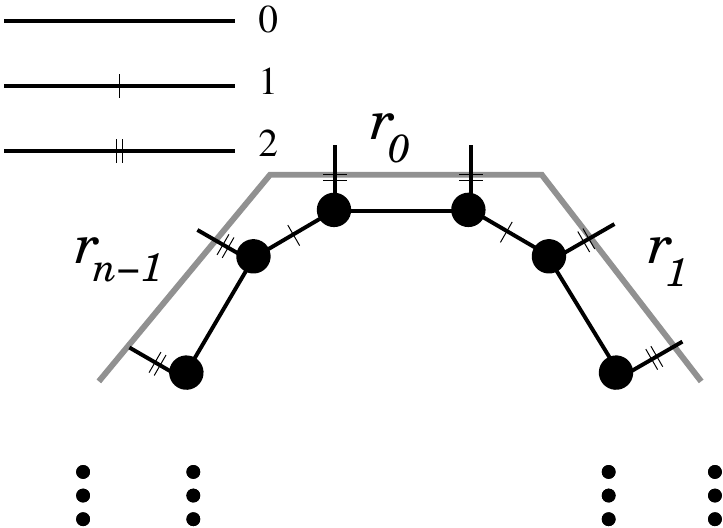}
  \caption{The voltage operator $(\cY,W^n,\eta_1)$.}
  \label{fig:fam1op2}
\end{figure}

Using \cref{teo:thetacomp}, we notice that we can interpret
$\cO$ as the result of concatenating a different operation $\cO'$ with
the opposite operation. The operation $\cO'$ is given by the voltage
operator $(\cY,W^n,\eta_1)$, obtained by replacing the edges of color
2 in $\cX$ by semi-edges with the same voltage (see
\cref{fig:fam1op2}).

If we want to realize an $n$-maniplex as a topological object, what we
usually do is interpret each flag $\Phi$ as an $(n-1)$-simplex. We
label the vertices of this simplex with the numbers from $0$ to $n-1$
and we label each facet of the simplex according to its opposite
vertex. Finally, for every two $i$-adjacent flags, we identify their
facets with the label $i$ in a way that vertices with the same label
are also identified. If $F$ is an $i$-face, for each flag $\Phi$ in
$F$ we consider the $i$-simplex defined by the vertices of the
$(n-1)$-simplex associated to $\Phi$ whose label is smaller or equal
to $i$, and we interpret $F$ as the union of such $i$-simplices. By
doing this, the vertex labeled $i$ of the simplex corresponding to the
flag $\Phi$ corresponds to a fixed point in the interior of the
$i$-face $(\Phi)_i$, which is the same for every other flag in that
same $i$-face. We call this point the \emph{barycenter of $(\Phi)_i$}.

\cref{fig:fam1op2} is telling us that to apply the operation $\cO'$ we
have to replace each flag by an $n$-gon, and that $i$-adjacent flags
correspond to $n$-gons sharing their $i$-th side. One way to do this
topologically is as follows: We start by representing each flag as a
simplex, as usual. Then, for each $i\in \{0,1,\ldots,n-1\}$ we draw a
point $u_i$ at the barycenter of the intersection of the facets
numbered $i$ and $i-1$ (reduced modulo $n$ if necessary). Then we
construct an $n$-gon by drawing a straight segment between $u_i$ and
$u_{i+1}$ for each $i$. Notice that the edge joining $u_i$ to
$u_{i+1}$ is contained in the facet with the label $i$. Therefore, the
polygons corresponding to a pair of $i$-adjacent flags will share
their $i$-th sides, just as requested.

\begin{figure}[h!]
  \centering
  \includegraphics[width=0.5\linewidth]{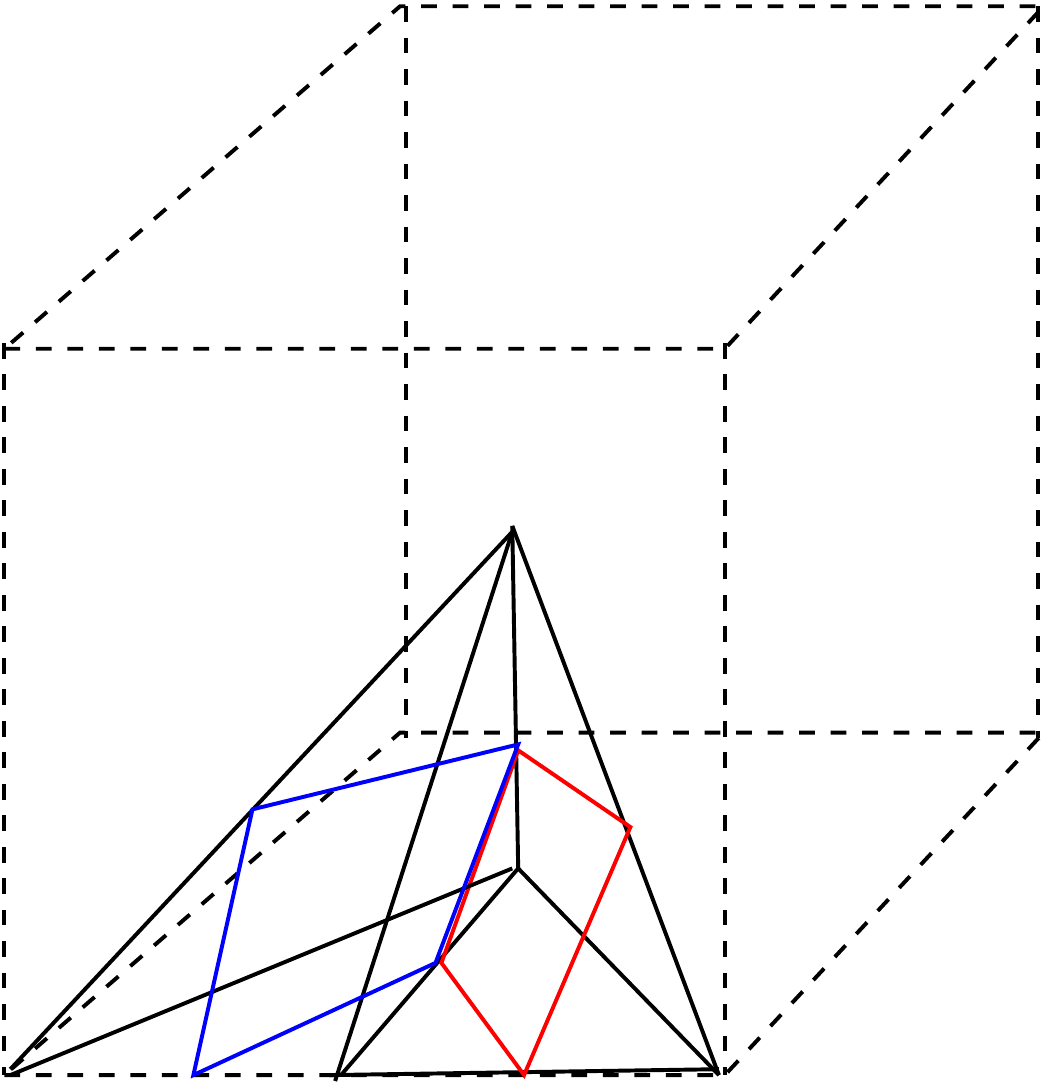}
  \caption{The operation $\cO'$ in rank 4.}
  \label{fig:opgeom1}
\end{figure}

For $n=3$, the intersection of the facet $i$ and the facet $i-1$ is
exactly the vertex $i+1$, so the triangle we get by applying $\cO'$ is
the same as the triangle representing the flag. In other words, in
rank 3, the operation $\cO'$ is the barycentric subdivision. In rank 4
the polygons we get have a vertex in the midpoint of each edge of a
particular Petrie polygon of the simplex. In \cref{fig:opgeom1} we see
an example in rank 4. Here the facets of our 4-polytope are cubes, and
we illustrate the simplices corresponding to two 0-adjacent flags
together with the polygons that they turn into when we apply the
operation $\cO'$.

Note that these operations can be applied to any premaniplex of the
correct rank, but we can only be certain that the result has symmetry
type in family 1 when the operated object is a regular polyhedron that
is non-degenerate, non-self-dual, and has simple Petrie polygons.

\subsection{Family 2}\label{sec:VoltOps2}

\begin{figure}[h!]
  \centering
  \includegraphics[width=0.5\linewidth]{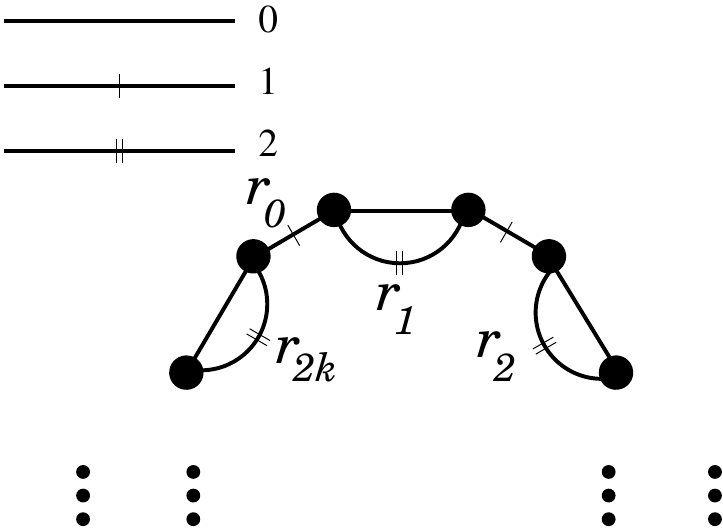}
  \caption{The voltage operator  for family 2.}
  \label{fig:fam2op}
\end{figure}

Now let us consider the construction of \cref{teo:fam2const}. Remember
that this construction starts with a chiral polytope $\cP$ of rank
$2k+1$ and gives us a polyhedron $\cX^\xi$ with $4k$-gonal facets.
The polytope $\cP$ can be constructed as $\cZ^\nu$ where $\cZ$ is the
$(2k+1)$-premaniplex with two vertices $a$ and $b$ joined by an edge
of each color from $0$ to $2k$, and the voltage assignment $\nu$ is
given by $\nu(^ia)=\tau_i$. Let $\cO$ be the voltage operator in
\cref{fig:fam2op}. It is similar to the operation in
\cref{fig:fam1op}, but the voltages of the edges of color 2 go from
$r_1$ to $r_k$ and the edge of color 1 between the first and last
flags in the cycle, has voltage $r_0$. We have removed the polygon in
gray to keep the image cleaner.

Using \cref{teo:theta} we see that the maniplex $\cX^\xi$ constructed
in \cref{teo:fam2const} is isomorphic to the result of applying the
operation $\cO$ to the chiral polytope $\cP\cong\cZ^\nu$.

\begin{figure}[h!]
  \centering
  \includegraphics[width=0.5\linewidth]{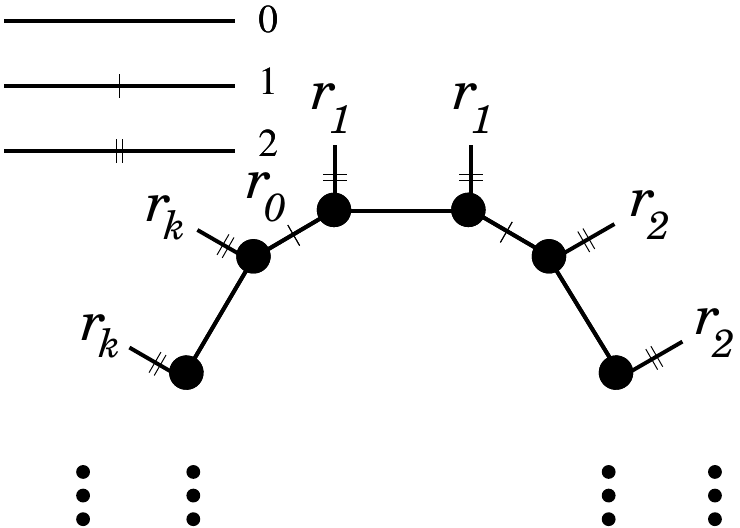}
  \caption{The operator for $\cO'$.}
  \label{fig:fam2op2}
\end{figure}

Again, the operation $\cO$ can be interpreted as the concatenation of
another operation $\cO'$ and the opposite operation, which turns the
edges of color 2 into semi-edges (see \cref{fig:fam2op2}).

\begin{figure}[h!]
  \centering
  \includegraphics[width=0.5\linewidth]{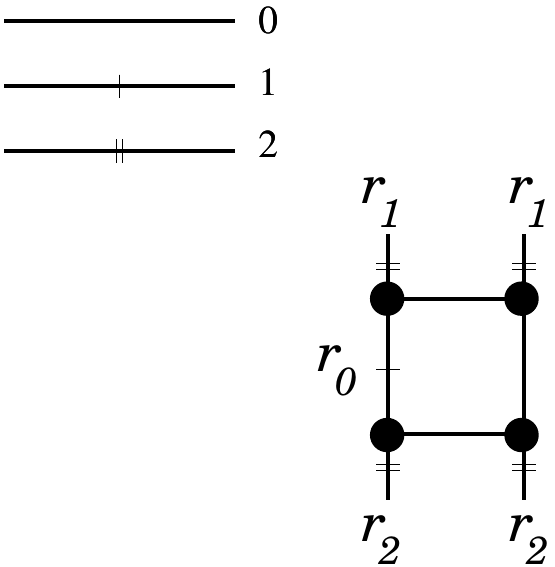}
  \caption{The operator for $\cO'$ in rank 3.}
  \label{fig:fam2op3}
\end{figure}

Now we would like to give a topological interpretation for the
operation $\cO'$, however this is hard to do in high rank, so we will
only do so for the case when $\cO'$ turns a 3-polytope (a
polyhedron) into a polyhedron with square faces. The corresponding
voltage operator appears in \cref{fig:fam2op3}. Even though this
particular operation does not give examples with trivial facet
stabilizer (as the vertical reflection in~\cref{fig:fam2volts} lifts),
understanding this operation in rank 3 might be the first step towards
understanding the ones in higher rank.

One can check that the operation $\cO'$ in rank 3 turns every pair of
0-adjacent flags in $\cP$ into a square. Such a pair determines an
edge (1-face) and a facet of $\cP$. Pairs of 0-adjacent flags that
share a 1-face are turned into squares that share two opposite sides,
while pairs that have 1-adjacent flags turn into squares that share 1
side.

A topological way to do this is as follows: Replace each edge $e$ of
$\cP$ with a cylinder split in half by a plane containing its
axis. This divides the cylinder into two squares (making the ``walls''
of a digonal prism). We can label each square with one of the two
facets of $\cP$ containing $e$. One can also label each end of the
cylinder with one of the vertices incident to $e$. Then, given a flag
$\Phi$ of $\cP$ there is a semicircle that corresponds to taking the
cylinder corresponding to $(\Phi)_1$, taking the half labeled by
$(\Phi)_2$, and taking the side of that square that corresponds to
$(\Phi)_0$. Then for each flag $\Phi$ we glue the semicircle
corresponding to $\Phi$ to the one corresponding to $r_1 \Phi$. By
doing this we get a map on the surface obtained by replacing the
1-faces of $\cP$ by cylinders. This resulting map is the polyhedron
$\cO'(\cP)$.

\begin{figure}[h!]
  \centering
  \includegraphics[width=0.7\textwidth]{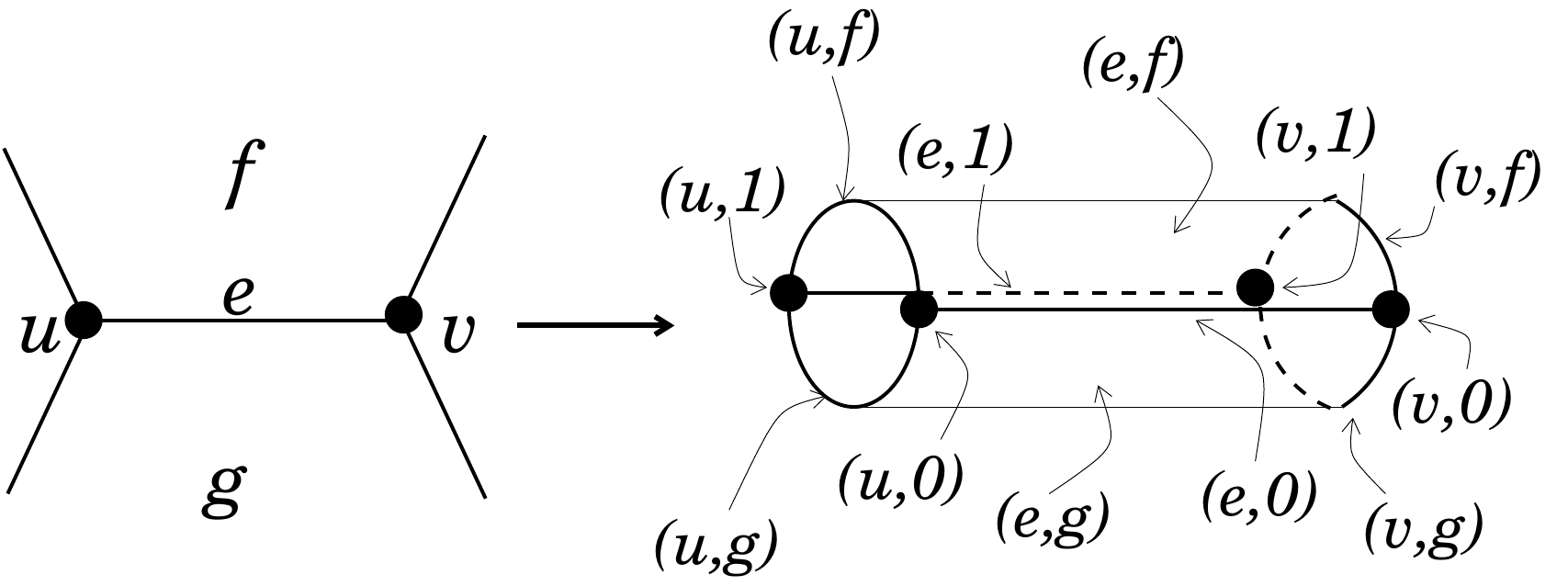}
  \caption{The operation $\cO'$ in rank 3.}
  \label{fig:opgeom2}
\end{figure}

We can interpret this operation also in the context of incidence
geometries. Recall that $(\cP)_i$ denotes the set of $i$-faces of
$\cP$. The vertices of $\cO'(\cP)$ can be thought of as
$(\cP)_0\times \{0,1\}$. The edges of $\cO'(\cP)$ can be thought of as
$(\cP)_1\times \{0,1\}$, together with the pairs of the form $(v,f)$
where $v$ is a vertex of $\cP$ and $f$ is a facet of $\cP$ incident to
$v$. The edge $(e,\varepsilon)$ is incident to the vertex $(v,\delta)$
if and only if $\varepsilon=\delta$ and $e$ is incident to $v$ in
$\cP$, while the edge $(v,f)$ is incident to $(v,0)$ and $(v,1)$. The
facets of $\cO'(\cP)$ correspond to pairs $(e,f)$ where $e$ is an edge
of $\cP$ and $f$ is a facet of $\cP$ incident to $e$. The facet
$(e,f)$ is incident to the edges $(e,0)$, $(e,1)$, $(v,f)$ and
$(u,f)$, where $v$ and $u$ are the vertices of $\cP$ incident to $e$
(see \cref{fig:opgeom2}).

One can verify that the voltage operator corresponding to this
operation is exactly the one in \cref{fig:fam2op3}.

With this interpretation it is also clear that there will always be an
automorphism of $\cO'(\cP)$ that interchanges $(v,0)$ with $(v,1)$ for
every vertex $v$, interchanges $(e,0)$ with $(e,1)$ for every edge $e$
and fixes every other element. In particular, the facet stabilizer in
$\cO'(\cP)$ (and in $\cO(\cP)$) is not trivial.


\section{Alternating semiregular polyhedra with trivial facet
  stabilizer}\label{sec:AltSemiReg}

\begin{figure}[h!]
  \centering
  \includegraphics[width=0.5\linewidth]{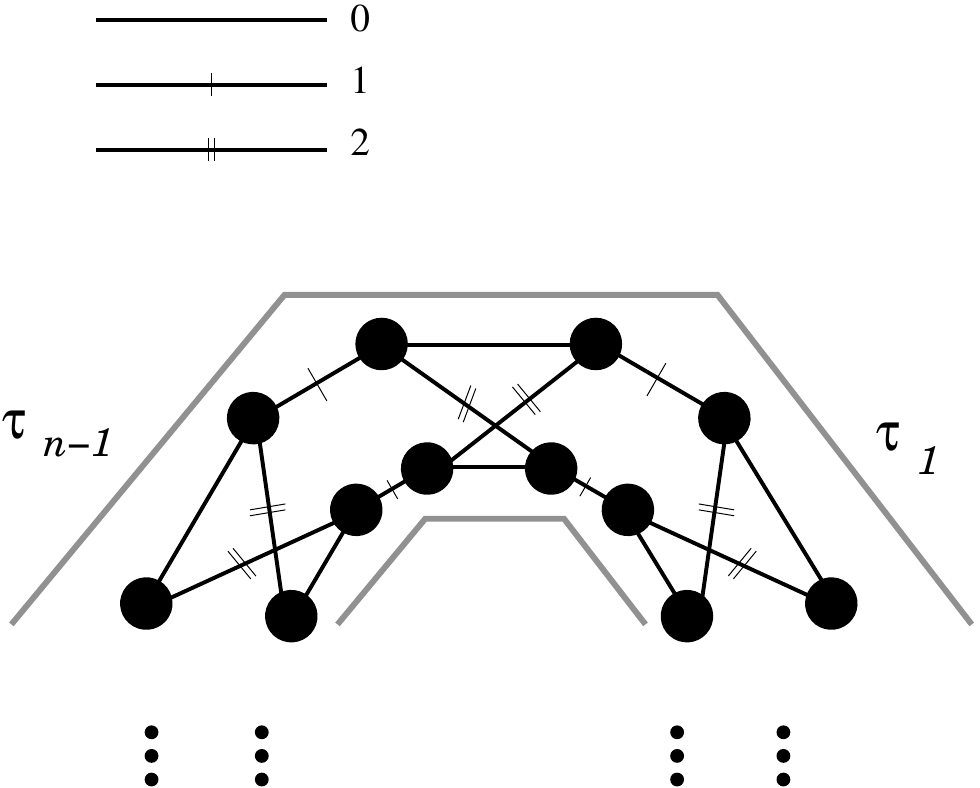}
  \caption{The voltage graph resulting from applying the operation of
    \cref{sec:VoltOps1} to $2^n_\emptyset$.}
  \label{fig:AltSemiReg}
\end{figure}

In \cref{sec:VoltOps1} we have discussed how polyhedra with symmetry
type graph in family 1 can be obtained by applying the voltage
operation illustrated in \cref{fig:fam1op} to a regular polytope of
rank $n$. If instead we apply this same operation to a chiral polytope
$\cP$ with automorphism group
$\Gamma(\cP)=\gen{\tau_1,\tau_2,\ldots,\tau_{n-1}}$, where once again
$\Phi\tau_i = r_i r_0\Phi$ for all $i$ and some fixed base flag
$\Phi$, then \cref{teo:theta} tells us that the result is the derived
graph of the voltage graph in \cref{fig:AltSemiReg}. This graph is
obtained by taking two $n$-gons and then connecting corresponding
edges with a ``cross'' of edges of color 2. That is, the left flag of
the $i$-th edge of one copy is 2-adjacent to the right flag of the
$i$-th edge of the second copy and vice-versa. For $i>0$ the darts of
color $2$ that start on the $i$-th edge of the first copy have voltage
$\tau_i$, and the voltage of the darts of color 2 starting at the 0-th
edge of any copy have trivial voltage. The way to construct this graph
is to apply the voltage operation in question to the symmetry type
graph of a chiral $n$-polytope, which is $2^n_\emptyset$. See
\cref{teo:theta} to see how to determine the voltage assignment.

Let $(\cX,\Gamma(\cP),\xi)$ be the voltage graph in
\cref{fig:AltSemiReg}. Observing that only the edges of color two
can have non-trivial voltage we can see that the facets of the derived
graph will all be $n$-gons. If we remove the edges of color 2 from
$\cX$ we get two connected components. This means that
$\Gamma(\cP)$ acts with two orbits on the set of facets. On the other
hand, the number of connected components after removing the edges of
color 0 depends on the parity of $n$: for odd $n$ we get one connected
component, while for even $n$ we get two connected components.

Notice that, since every edge of color 2 connects the two
components of $\cX_{\ol{2}}$, we know that every 1-face (edge) of
$\cX^\xi$ is shared by facets in different orbits.

Let $\phi$ (lowercase) be the right flag of the edge $e_0$ of the
outer $n$-gon. We will use $\phi$ as the base flag in $\cX$. One can
find a closed path $W_i$ based at $\phi$ with voltage $\tau_i$ as
follows: First go from $x$ to a flag in the edge $e_i$ of the outer
$n$-gon by using darts of colors 0 and 1. Second, use a dart $d_i$ of
color 2 to get to the edge $e_i$ of the inner $n$-gon. Third, use
edges of colors 0 and 1 to get to the left flag of the edge $e_0$ of
the inner $n$-gon. And fourth, use the dart $d_0$ of color 2 starting
at this flag to go back to $\phi$.

Notice that if $\sigma$ is an automorphism of $\cX$, then
$\xi(W_i\sigma) = \xi(d_0\sigma)\xi(d_i\sigma)$. 
Notice also that every closed path based at $\Phi$ can be factored
into paths of the form $W_i^{\pm 1}$. So, by knowing the voltage of
the paths $W_i$ we can calculate the voltage of any other closed path
based at $\phi$. Then, by \cref{teo:lifts}, $\sigma$ lifts if and only
if the mapping
$\sigma^*:\tau_i\mapsto \xi(d_0\sigma)\xi(d_i\sigma)$
extends to an automorphism of $\Gamma$.

We want to make sure that if $\sigma$ is a non-trivial automorphism of
$\cX$, then it does not lift. To do this, we first notice that the
involution $\tau_{n-1}$ is mapped by $\sigma^*$ to
$\xi(d_0\sigma)\xi(d_{(n-1)}\sigma)$. Then
$\tau_{n-1}\sigma^* = \tau_i^{\pm1}\tau_{i\pm 1}^{\pm 1}$ for some
$i\in\{0,\ldots,n-1\}$ and some choice of signs $\pm$ (subscripts are taken
modulo $n$). This product is an involution if and only if one of the
factors is $\tau_0$ (the identity) and the other one $\tau_{n-1}$
(this is again a consequence of the fact that there are no degenerate
chiral polytopes). Since $\sigma$ is not the identity, this leaves us
with three possibilities:

{\bf Case 1: $\sigma$ switches the inner and the outer $n$-gons
  (mapping $\phi$ to $\phi^{0,2}$).} In this case $\tau_i\sigma^* =
\tau_i^{-1}$. In this case, $\sigma^*$ cannot be extended to an
automorphism of $\Gamma$.

{\bf Case 2: $\sigma$ is the reflection on a line through the left
  corner of $e_0$ (mapping $\phi$ to $\phi^{0,1,0}$).} In this case
$\sigma^*$ maps $\tau_i$ to $\tau_{n-1}\tau_{n-1-i}$. We can see that
$\tau_{n-1}\tau_{n-1-i}$ maps $\Phi$ (the base flag of $\cP$) to
$\Phi^{0,n-1-i,n-1,0}$. This means that in the dual $\cP^*$ it maps
$\Phi$ to $\Phi^{n-1,i,0,n-1}$. If we now look at the action on
$\Phi^{n-1,0}$ in $\cP^*$ we see that
$\Phi^{n-1,0}(\tau_i\sigma^*) = \Phi^{n-1,i} =
(\Phi^{n-1,0})^{0,i}$. In other words, $\sigma^*$ maps each
distinguished generator of $\cP$ to the corresponding generator of
$\cP^*$ when we use $\Phi^{n-1,0}$ as the base flag. Then $\sigma^*$
can be extended to an automorphism of $\Gamma$, if and only if
$\cP$ is \emph{properly self-dual} (that is, there is an isomorphism
from $\cP$ to $\cP^*$ that maps each flag to a flag in the same
orbit, see~\cite[Section 2.2.2]{ACP}).

{\bf Case 3: $\sigma$ is the composition of the previous two cases
  (mapping $\phi$ to $\phi^{0,1,2}$).} In this case
$\tau_i\sigma^* = \tau_{n-1}\tau_{n-1-i}^{-1}$, which maps $\Phi$ to
$\Phi^{n-1-i,n-1}$. This means that in the dual $\cP^*$ it maps $\Phi$
to $\Phi^{i,0}$. If we now look at the action on $\Phi^0$ in $\cP^*$
(which is $\Phi^{n-1}$ in $\cP$) we see that
$\Phi^0(\tau_i\sigma^*) = \Phi^i = (\Phi^0)^{0,i}$. In other words,
$\sigma^*$ maps each distinguished generator of $\cP$ to the
corresponding generator of $\cP^*$ when we use $\Phi^0$ as the base
flag. Then, $\sigma^*$ could be extended to an automorphism of
$\Gamma$, if and only if $\cP$ is \emph{improperly self-dual} (that
is, there is an isomorphism from $\cP$ to $\cP^*$ that maps each flag
to a flag in the opposite orbit, see~\cite[Section 2.2.2]{ACP}).

Therefore we get the following lemma:

\begin{lema}\label{lema:quiral-auto-dual}
  The maniplex $\cX^\xi$ derived from the voltage graph $(\cX,\xi)$ in
  \cref{fig:AltSemiReg} has symmetry type graph $\cX$ if and only if
  $\cP$ is not self-dual.
\end{lema}

In particular, if $n$ is odd then $\cX^\xi$ is an alternating
semiregular maniplex, and if $\cP$ is not self-dual, then $\cX^\xi$
has trivial facet stabilizer. Now we want to prove that it is also
polytopal.

To do this we will use the following generalization of
\cref{teo:PolyPreExt} (the proof is identical):

\begin{teo}\label{teo:PolyPreExtGen}
  Let $(\cX,\Gamma,\xi)$ be a voltage $(r+1)$-premaniplex. Suppose
  $\xi(d)=1$ for every dart $d$ of color smaller than $r$. Then
  $\cX^\xi$ is a polytope if and only if the connected components of
  $\cX_{\ol{r}}$ are polytopes, and whenever there is a path $W$ with
  voltage 1 and colors in $[m,r]$ ending in the same connected
  component of $\cX_{\ol{r}}$ as it started, there is also a path $W'$
  with the same endpoints as $W$ that uses only colors in $[m,r-1]$.
\end{teo}

If we start a path $W$ at $\phi$ (the base flag of $\cX$) and follow
darts that alternate between colors 1 and 2 (starting with color 1),
we will see that its voltage is $\prod_{i=1}^j \tau_i^{\epsilon_i}$,
where the subscript of $\tau_i$ is taken modulo $n-1$ (using
$\tau_{n-1}$ when $i\equiv 0 \mod n-1$), $\epsilon_i$ is $-1$ whenever
$i\equiv n+1 \mod 2n$ and 1 in any other case and $j$ is the number of
edges of color 2 with non-trivial voltage used by $W$. By using
\cref{prop:prod_tau} we get that the voltage of $W$ acts like a
$2t$-step rotation in the Petrie polygon of $\Phi$ when
$j=2(n-1)t$. Using the same argument as in \cref{teo:fam2const} we get
that if $W$ has trivial voltage, then $j$ is divisible by
$2(n-1)$. This means that the last dart with non-trivial voltage used
by $W$ is the dart of color 2 ending on the right flag of the
$(n-1)$-th edge of the inner $n$-gon, which is $r_1
r_2\phi$. Therefore, $W$ must end in either $r_1 r_2\phi$, $r_2\phi$,
$\phi$ or $r_1\phi$, but only the last two are in the same connected
component of $\cX_{\ol{2}}$ as $\phi$. So, by using
\cref{teo:PolyPreExtGen} we get the following result.

\begin{teo}
  Let $n$ be an odd natural number and $\cP$ be a non-self-dual chiral
  $n$-polytope with simple Petrie polygons. Let $\cO$ be the voltage
  operation of \cref{sec:VoltOps1}. Then $\cO(\cP)$ is an alternating
  semiregular polyhedron with trivial facet stabilizer and its
  symmetry type graph is the one in \cref{fig:AltSemiReg}.
\end{teo}

\section{Higher rank}\label{sec:highrank}

In~\cite{OrbitasEnSemiregs}, given an $n$-polytope $\cP$, the authors
construct a vertex-transitive polytope $2^{\cP}$ whose vertex-figures
are isomorphic to $\cP$. If all the facets of $\cP$ are isomorphic to
a regular $(n-1)$-polytope $\cK$, then $2^{\cP}$ is semiregular and
all its facets are isomorphic to $2^{\cK}$. If $\cP$ has trivial
automorphism group, then $2^\cP$ has trivial vertex stabilizer, which
implies that $2^\cP$ has as many flag-orbits as $\cP$ has flags. The
authors show that for $n\geq 3$ there are arbitrarily large polytopes
with trivial automorphism group and isomorphic regular facets,
therefore showing that there are semiregular $(n+1)$-polytopes with
arbitrarily many flag-orbits. However, $2^\cP$ has also as many
facet-orbits as $\cP$ has facets.

Using \cref{teo:fam1const} we can construct examples of 3-polytopes
with trivial facet stabilizer (therefore with as many flag orbits as
flags in a facet) that are also facet-transitive. In this section we
will generalize family 1 to construct some examples of semiregular
maniplexes of higher rank with trivial facet stabilizer and exactly
one orbit on facets.

Let $\cM$ be a regular $n$-maniplex with facets isomorphic to an
$(n-1)$-maniplex $\cK$. Suppose that $\Gamma(\cK)$ has an involution
$\sigma$ in its center. Given a flag $\Phi$ in a facet $F$ of $\cM$ we
define $r_n \Phi$ as $\Phi\tau_F^{-1} \sigma \tau_F$, where $\tau_F$
is some isomorphism from $\cK$ to $F$. One can prove that $r_n$ does
not depend on the choice of $\tau_F$ (however, it does depend on the
choice of $\sigma$). We may call $r_n\Phi$ the \emph{opposite flag of
  $\Phi$ in the facet $F$}. Moreover, if $G$ is an $i$-face incident
to $F$ with $i<n-1$, we can define the \emph{opposite face of $G$ in
  $F$} in a natural way. By adding edges of color $n$ between each
flag and its opposite (in its respective facet), we get an
$(n+1)$-premaniplex $\cX$ (in fact, it is usually a maniplex). If
$\cX_{\ol{0}}$ is connected, then we can use $\cX$ to construct a
semiregular maniplex with symmetry type $\cX$:

Let $F_0, F_1,\ldots, F_{N-1}$ be the facets of $\cM$. Let $\cP$ be a
regular polytope of rank $N$. Let
$\Gamma(\cP)=\gen{\rho_0,\ldots,\rho_{N-1}}$. We construct a voltage
assignment $\xi:D(\cX)\to \Gamma(\cP)$ by giving trivial voltage to
all darts of color smaller than $n$, and giving voltage $\rho_i$ to
the darts of color $n$ in the facet $F_i$. Then $\cX^\xi$ will be an
$(n+1)$-maniplex whose facets are isomorphic to $\cM$. If no
automorphism of $\cX$ induces an automorphism of $\Gamma(\cP)$, we can
be sure that $\cX$ is the symmetry type graph of $\cX^\xi$, implying
that $\cX^\xi$ has trivial facet stabilizer. If $\cP$ is
non-degenerate, the only possible non-trivial automorphism induced by
a permutation of the generators of $\cP$ would be one that maps
$\rho_i$ to $\rho_{N-1-i}$. In most cases, we will be able to find a
labeling of the facets of $\cM$ in a way that there is no automorphism
that induces this permutation. In the rare case that this happens to
be impossible, by demanding $\cP$ to not be self-dual we make sure
that this permutation of generators does not induce an automorphism of
$\Gamma(\cP)$.

Just as we did in \cref{sec:VoltOps1}, we may interpret this
construction as applying some voltage operation that turns regular
$N$-polytopes into semiregular, facet-transitive $(n+1)$-maniplexes
with trivial facet stabilizer. Applying this operation to a chiral
$N$-polytope would give an alternating semiregular $(n+1)$-maniplex
with trivial facet stabilizer (just as in \cref{sec:AltSemiReg}).

The problem comes with polytopality. We can be certain that we need
$\cM$ to be a polytope for $\cX^\xi$ to be a polytope. However it is
not obvious what other conditions are necessary or sufficient to make
sure that $\cX^\xi$ is a polytope. This is an interesting problem to
explore in the future. We have tried doing some experimentation with
RAMP, but the examples become too big too quickly.

\section*{Acknowledgments}
This paper was completed while the author had a postdoctoral position
at the Department of Mathematics of Northeastern University. The
author wishes to thank Gabe Cunningham for his help with RAMP. The
author also thanks the referees for their observations, in particular
referee Y for helping on finding a mistake on a previous version of
\cref{teo:fam2const}.

\bibliographystyle{plain}
\bibliography{Bibliografia}

\end{document}